%&amstex
% This is an AMS-TeX file and should be compiled
% using AMS-TeX.
% the command should be:  amstex file  or something
% like that, depending on what has been set up
% on the system you're using. 
\magnification=1200

\loadmsam
\loadmsbm
\loadeufm
\loadeusm
\UseAMSsymbols
%\hsize=6.00 true in
%\hoffset=.50 true in
%\voffset=-0.1 true in
%\vsize=8.75 true in

\font\BIGtitle=cmr10 scaled\magstep3
\font\bigtitle=cmr10 scaled\magstep1
\font\boldsectionfont=cmb10 scaled\magstep1
\font\section=cmsy10 scaled\magstep1

\def\scr#1{{\fam\eusmfam\relax#1}}

\def\scrB{{\scr B}}

\def\scrS{{\scr S}}

\def\gr#1{{\fam\eufmfam\relax#1}}

	\def\grg{{\gr g}}
	
\def\grI{{\gr I}}

\def\grL{{\gr L}}

\def\db#1{{\fam\msbfam\relax#1}}

\def\dbA{{\db A}} 
\def\dbC{{\db C}} 
 \def\dbF{{\db F}}
\def\dbG{{\db G}}

 \def\dbN{{\db N}}
 
\def\dbQ{{\db Q}}

 \def\dbZ{{\db Z}}

\def\Ebar{\bar{E}}

\def\kbar{\bar{k}}
\def\lbar{\bar{l}}

\def\Ker{\text{Ker}}

\def\der{\text{der}}

\def\sc{\text{sc}}
\def\Res{\text{Res}}
\def\ab{\text{ab}}
\def\ad{\text{ad}}

\def\Gal{\text{Gal}}
\def\Hom{\text{Hom}}
\def\End{\text{End}}

\def\Spec{\text{Spec}}

\def\Lie{\text{Lie}}

\def\leaderfill{\leaders\hbox to 1em
     {\hss.\hss}\hfill}
\def\nspace{\lineskip=1pt\baselineskip=12pt\lineskiplimit=0pt}

     %the way to use this is "\Proclaim{Theorem 1.1.}" for instance.
\def\finishproclaim{\par\rm
     \ifdim\lastskip<\medskipamount\removelastskip
     \penalty55\medskip\fi}
\def\proof{\par\noindent {\it Proof:}\enspace}

\def\Ref[#1]{\par\hang\indent\llap{\hbox to\parindent
     {[#1]\hfil\enspace}}\ignorespaces}
\def\Item#1{\par\smallskip\hang\indent\llap{\hbox to\parindent
     {#1\hfill$\,\,$}}\ignorespaces}
\def\ItemItem#1{\par\indent\hangindent2\parindent
     \hbox to \parindent{#1\hfill\enspace}\ignorespaces}

\def\Le{{\mathchoice{\,{\scriptstyle\le}\,}
  {\,{\scriptstyle\le}\,}
  {\,{\scriptscriptstyle\le}\,}{\,{\scriptscriptstyle\le}\,}}}
\def\Ge{{\mathchoice{\,{\scriptstyle\ge}\,}
  {\,{\scriptstyle\ge}\,}
  {\,{\scriptscriptstyle\ge}\,}{\,{\scriptscriptstyle\ge}\,}}}

\def\arrowsim{\,\smash{\mathop{\to}\limits^{\lower1.5pt
  \hbox{$\scriptstyle\sim$}}}\,}

\def\doublemaprights#1#2#3#4{\raise3pt\hbox{$\mathop{\,\,\hbox to
     #1pt{\rightarrowfill}\kern-30pt\lower3.95pt\hbox to
     #2pt{\rightarrowfill}\,\,}\limits_{#3}^{#4}$}}

\def\rightcapdownarrow{\raise9pt\hbox{$\ssize\cap$}\kern-7.75pt
     \Big\downarrow}

\def\rcapmapdown#1{\rightcapdownarrow\kern-1.0pt\vcenter{
     \hbox{$\scriptstyle#1$}}}

\def\rmapdown#1{\Big\downarrow\kern-1.0pt\vcenter{
     \hbox{$\scriptstyle#1$}}}
\def\rightsubsetarrow#1{{\ssize\subset}\kern-4.5pt\lowe r2.85pt
     \hbox to #1pt{\rightarrowfill}}
\def\longtwoheadedrightarrow#1{\raise2.2pt\hbox to #1pt{\hrulefill}
     \!\!\!\twoheadrightarrow}

\def\Gal{\operatorname{\hbox{Gal}}}
\def\Hom{\operatorname{\hbox{Hom}}}

\def\im{\hbox{Im}}

\NoBlackBoxes
\parindent=25pt
\document
\footline={\hfil}

\null
\noindent
manuscripta math. 112, 325--355 (2003) $\;\;\;\;\;\;\;\;\;\;\;\;\;\;\;\;\;\;\;\;\;\;\;\;\;\;\;\;\;\;\;\;\;\;\;\;\;\;\;\;\;\;\;$ Springer-Verlag 2003
\bigskip
\noindent
{\bigtitle Adrian Vasiu}
\vskip 0.15in
\noindent
{\BIGtitle Surjectivity Criteria for p-adic Representations, Part I}
\footline={\hfill}
\vskip 0.15in
\noindent
{Received: 10.09.2002}
\noindent
{Revised version: 27.01.2003}
\noindent
{Manuscripta Mathematica, MM-No:577}
\smallskip\noindent
Published online: 15.10.2003
\vskip 0.2in
\noindent
{\bf ABSTRACT}. We prove general surjectivity criteria for $p$-adic representations. In particular, we classify all adjoint and simply connected group schemes $G$ over the Witt ring $W(k)$ of a finite field $k$ such that the reduction epimorphism $G(W_2(k))\twoheadrightarrow G(k)$ has a section.
\noindent
$\vfootnote{} {Adrian Vasiu, Mathematics Department, University of Arizona, 617 N. Santa Rita, P.O. Box 210089, Tucson, AZ-85721, USA. e-mail: adrian\@math.arizona.edu}$
\noindent
$\vfootnote{} {{\it Mathematics Subject Classification (2000)}: Primary 11S23, 14L17, 17B45, 20G05 and 20G40}$

\vskip0.1in

\footline={\hss\tenrm \folio\hss}
\pageno=1

\bigskip
\noindent
{\boldsectionfont \S1. Introduction}

\bigskip
Let $p$ be a rational prime. Let $d$, $r\in\dbN$. Let $q:=p^r$. Let $A$ be an abelian variety of dimension $d$ over a number field $E$. The action of $\Gal(\Ebar/E)$ on the $\Ebar$-valued points of $A$ of $p$-power order defines a $p$-adic representation 
$$\rho_{A,p}:\Gal(\Ebar/E)\to GL(T_p(A))(\dbZ_p)=GL_{2d}(\dbZ_p),$$ 
where $T_p(A)$ is the Tate module of $A$. If $A$ is a semistable elliptic curve over $\dbQ$ (so $d=1$ and $E=\dbQ$), then Serre proved that for $p\Ge 3$ the homomorphism $\rho_{A,p}$ is surjective iff its reduction mod $p$ is irreducible (cf. [23, Prop. 21], [24, p. IV-23-24] and [26, p. 519]). 

Similar results are expected to hold for many homomorphisms of the form
$$\rho:\Gal(E)\to G(W(k)),$$
with $G$ a reductive group scheme over the Witt ring $W(k)$ of $k:=\dbF_{q}$ (like when the Zariski closure $G_{A,p}$ of $\im(\rho_{A,p})$ in $GL(T_p(A))$ is a reductive group scheme and $\rho$ is the factorization of $\rho_{A,p}$ through $G_{A,p}(\dbZ_p)$).  Cases involving $q=p$ and groups $G$ which are extensions of $\dbG_m$ by products of Weil restrictions of $SL_2$ groups were treated for $p\Ge 5$ in [22]. If $f:G\hookrightarrow GL_n$ is a monomorphism of interest, then in general the assumption that the reduction mod $p$ of $f(W(k))\circ\rho$ is an irreducible representation over $k$, is too weak to produce significant results. So the goal of this paper is to get criteria under which $\rho$ is surjective if its reduction mod $p$ (or occasionally mod $4$ for $p=2$) is surjective. 

A reductive group scheme $F$ over a connected affine scheme $\Spec(R)$ is assumed to have connected fibres. Let $F^{\der}$, $Z(F)$, $F^{\ab}$ and $F^{\ad}$ be the derived group, the center, the maximal commutative quotient and respectively the adjoint group of $F$. So $F^{\ad}=F/Z(F)$ and $F/F^{\der}=F^{\ab}$. Let $Z^0(F)$ be the maximal torus of $Z(F)$. Let $F^{\sc}$ be the simply connected semisimple group cover of $F^{\der}$. Let $c(F^{\der})$ be the degree of the central isogeny $F^{\sc}\to F^{\der}$. If $S$ is a closed subgroup of $F$, then $\Lie(S)$ is the $R$-Lie algebra of $S$. We now review the contents of this paper. 

\medskip\noindent
{\bf 1.1. Adjoint representations.} If the reduction mod $p$ of $\rho$ is surjective, then the study of the surjectivity of $\rho$ is intimately interrelated with the study of the adjoint representation $AD_{G_k}$ of $G(k)$ on $\Lie(G_k)$. Though not stated explicitly, this principle is present in disguise in [24, p. IV-23-24] and [21, 2.1]. In order to apply it in \S4, in \S3 we assume $G$ is semisimple, we work just with $G_k$ (without mentioning $G$), and we deal with the classification of subrepresentations of $AD_{G_k}$. If $G_k$ is split, simply connected and $G^{\ad}_k$ is absolutely simple, then such a classification was obtained in [13]. The main goal of \S3 is to extend loc. cit. to the more general context of Weil restrictions of scalars (see [2, 7.6]) of semisimple groups having absolutely simple adjoints (see 3.4, 3.10 and 3.11). The methods we use are similar to the ones of [13] except that we rely more on the work of Curtis and Steinberg on representations over $\bar k$ of finite groups of Lie type and on the work of Humphreys and Hogeweij on ideals of $\Lie(G^{\der}_{\bar k})$ (see [17] and [14]; see also [18, 0.13] and [20, \S 1]). Whenever possible we rely also on [13]. Basic properties of Lie algebras and Weil restrictions are recalled in 2.2 and 2.3. 

\medskip\noindent
{\bf 1.2. The problem.} For $s\in\dbN$ let $W_s(k):=W(k)/p^sW(k)$. Let $\Lie_{\dbF_p}(G_k)$ be $\Lie(G_k)$ but viewed just as an abelian group identified with $\Ker(G(W_{2}(k))\to G(k))$. So if $H$ is a normal subgroup of $G_k$, then we also view $\Lie_{\dbF_p}(H)$ as a subgroup of $G(W_2(k))$. Let $K$ be a closed subgroup of $G(W(k))$ surjecting onto $(G/Z^0(G))(k)$ (we think of it as the image of some $\rho$). Let $K_2:=\im(K\to G(W_2(k)))$. The problem we deal with is to find conditions which imply that $K$ surjects onto $(G/Z^0(G))(W(k))$. It splits into two cases: $G$ is or is not semisimple. In this Part I we deal with the first case and in Part II we will deal with the second case and with applications of it to abelian varieties. The group $G$ is semisimple iff the torus $Z^0(G)$ is trivial. If $G$ is semisimple and $p|c(G)$, then there are proper, closed subgroups of $G(W(k))$ surjecting onto $G(k)$ (cf. 4.1.1). On the other hand, we also have the following general result: 

\medskip\noindent
{\bf 1.3. Main Theorem.} {\it We assume that $G$ is semisimple, that $g.c.d.(p,c(G))=1$ and that $K$ surjects onto $G(k)$. We also assume that one of the following five conditions holds:

\medskip
{\bf a)} $q\Ge 5$;

\smallskip
{\bf b)} $q=3$ and for each normal subgroup $H$ of $G_k$ which is a $PGL_2$ or an $SL_2$ group we have $\Lie_{\dbF_3}(H)\cap K_2\neq \{0\}$; 

\smallskip
{\bf c)} $q=4$ and $G^{\ad}_k$ has no simple factor which is a $PGL_2$ group;

\smallskip
{\bf d)} $q=2$ and the following two additional conditions hold:

-- no simple factor of $G^{\ad}_k$ is a $PGL_2$ group, a Weil restriction from $\dbF_4$ to $k$ of a $PGL_2$ group, or an $SU_4^{\ad}$ group;

-- for each normal subgroup $H$ of $G_k$ whose adjoint is a $PGL_3$ group, a $SU_3^{\ad}$ group or split of $G_2$ Lie type, we have $\Lie_{\dbF_2}(H)\cap K_2\neq\{0\}$;

\smallskip
{\bf e)} $p=2$ and $K_2=G(W_2(k))$.

\medskip
Then we have $K=G(W(k))$.}

\medskip
The case $G=SL_2$ for $p\Ge 3$ is due to Lenstra (unpublished computations with $3\times 3$ matrices). Serre proved the case $p\Ge 5$ for $SL_n$ and $Sp_{2n}$ groups over $\dbZ_p$ and the mod $8$ variant of e) for $SL_n$ groups over $\dbZ_2$ (see [24, IV] and [26, p. 52]). Most of the extra assumptions of b) to d) were known to be needed before (for instance, cf. [12, Sect. 4] in connection to d) for the $G_2$ Lie type). Though the case $q>>0$ of 1.3 is considered well known, we do not know any other concrete literature pertaining to 1.3.

\medskip\noindent
{\bf 1.4. On the proof of 1.3.} The proof of 1.3 is presented in 4.7. Its main ingredients are 4.3 to 4.5 and most of \S3. In 4.5 we list all isomorphism classes of adjoint groups $G=G^{\ad}$ for which the short exact sequence $0\to\Lie_{\dbF_p}(G_k)\to G(W_2(k))\to G(k)\to 0$ of abstract groups has a section, i.e. the epimorphism $G(W_2(k))\twoheadrightarrow G(k)$ has a right inverse. They are the eight ones showing up concretely in 1.3 b) to d). Though most of them are well known, we were not able to trace a reference pertaining to the complete classification of 4.5. The main idea of 4.5 is the following approach of inductive nature (see 4.4). Let $\gamma_G\in H^2(G(k),\Lie_{\dbF_p}(G_k))$ be the class defining the mentioned short exact sequence, with $\Lie_{\dbF_p}(G_k)$ viewed as a left $G(k)$-module via $AD_{G_k}$. If $\gamma_G=0$, then all images of restrictions of $\gamma_G$ are 0 classes and so $\gamma_{G_0^{\ad}}=0$ for any semisimple subgroup $G_0$ of $G$ normalized by a maximal torus of $G$ (see 4.3.4). But if there is a simple factor of $G$ whose isomorphism class is not in the list, then we can choose $G_0$ such that the direct computations of 4.4 show that $\gamma_{G_0^{\ad}}\neq 0$; so $\gamma_G\neq 0$. 

We now detail how 4.5 and \S3 get combined to prove 1.3 for $p\Ge 5$. Serre's method of [24, IV] can be adapted to get that for $p\Ge 5$ it is enough to show that $K_2=G(W_2(k))$ (see 4.1.2). Based on 4.5 we know that $\Lie_{\dbF_p}(G_k)\cap K_2$ is not included in $\Lie_{\dbF_p}(Z(G_k))$. So $\Lie_{\dbF_p}(G_k)\cap K_2=\Lie_{\dbF_p}(G_k)$, cf. 3.7.1 and 3.10 2). So $K_2=G(W_2(k))$. 

\medskip
Though \S 3 and \S 4  handle also the cases $p=2$ and $p=3$, \S 4 does not bring anything new to Serre's result recalled before 1.1; however, one can adapt our results to elliptic curves or to [22] in order to get meaningful results in mixed characteristics $(0,2)$ and $(0,3)$ (for instance, cf. 1.3 c) to e)). This Part I originated from seminar talks in Berkeley of Ribet and Lenstra; the first draft of \S 3 was a letter to Lenstra. 

\bigskip
\noindent                                             
{\boldsectionfont \S 2. Preliminaries}

\bigskip
In 2.1 we list our notations and conventions. In 2.2 and 2.3 we recall simple properties of Lie algebras and respectively of Weil restrictions of scalars.
 
\medskip\noindent
{\bf 2.1. Notations and conventions.} Always $n\in\dbN$. We denote by $k_1$ a finite field extension of $k=\dbF_{q}$. We abbreviate absolutely simple as a.s. and simply connected as s.c. Let $R$, $F$ and $S$ be as in \S 1. We say $F^{\ad}$ is simple (resp. is a.s.) if (resp. if each geometric fibre of) it has no proper, normal subgroup of positive relative dimension. If $M$ is a free $R$-module of finite rank, then $GL(M)$ ($SL(M)$, etc.) are viewed as reductive group schemes over $R$. So $GL(M)(R)$ is the group of $R$-linear automorphisms of $M$. 

Let $F$ be semisimple and $\Spec(R)$ connected. Let $o(F)$ be the order of $Z(F)$ as a finite, flat, group scheme. So $o(SL_2)=2$. We have $c(F)=o(F^{\sc})/o(F)$. See item (VIII) of [5, planches I to IX] for $o(F)$'s of s.c. semisimple groups having a.s. adjoints. If $X_*$ (resp. $X_{R_1}$ or $X$) is a scheme over $\Spec(R_1)$, then $X_{*R}$ (resp. $X_R$) is its pull back to $\Spec(R)$. All modules are left modules. A representation of a Lie algebra or a finite group $LG$ is also referred as an $LG$-module. See [25, p. 132] and [4, 16.3 to 16.6] for Lang theorem on connected, affine groups over finite fields. See [29] and [11] for the Lie types of semisimple groups over $k$ having a.s. adjoints. For Curtis and Steinberg theory of representations of finite groups of Lie type we refer to [9], [10] and [28, Theorem 1.3]; see also [11, 2.8] for an overview and [3, 6.4, 7.2 and 7.3] for the case of Chevalley groups. 

See [5, \S4] and [16, \S11] for the classification of connected Dynkin diagrams. We say $F$ is of (or has) isotypic $DT\in\{A_n,B_n,C_n|n\in\dbN\}\cup\{D_n|n\Ge 3\}\cup\{E_6,E_7,E_8,F_4,G_2\}$ Dynkin type if the connected Dynkin diagram of any simple factor of every geometric fibre of $F^{\ad}$ is $DT$; if $F^{\ad}$ is a.s. we drop the word isotypic. We use the standard notations for classical reductive group schemes over $k$, $W(k)$ or $\dbC$ (see [1]). So $PGL_n=GL_n^{\ad}=SL_n^{\ad}$, $PGSp_{2n}=Sp_{2n}^{\ad}=GSp_{2n}^{\ad}$, etc. Also $GU_n$ is the non-split form of $GL_n$ over $k$ or $W(k)$, $SU_n=GU_n^{\der}$ and $PGU_n:=SU_n^{\ad}=GU_n^{\ad}$. 

\medskip\noindent
{\bf 2.2. Lie algebras.} Let $x$ be an independent variable. As an $R$-module, we identify $\Lie(S)$ with the tangent space of $S$ at the identity section, i.e. with $\Ker(S(R[x]/x^2)\to S(R))$, where the $R$-epimorphism $R[x]/(x^2)\twoheadrightarrow R$ takes $x$ into $0$. If $y$, $z\in\Lie(S)$, then the Lie bracket $[y,z]$ is $yzy^{-1}z^{-1}$, the product being taken inside $\Ker(S(R[x]/x^2)\to S(R))$. We now assume that $S$ is a smooth group scheme over $R$. So $\Lie(S)$ is a free $R$-module of rank equal to the relative dimension of $S$. The representation of $S$ on $GL(\Lie(S))$ defined by inner conjugation is called the adjoint representation. Let $$AD_S:S(R)\to GL(\Lie_R(S))(R)$$ 
be the adjoint representation evaluated at $R$. 

\medskip\noindent
{\bf 2.3. Weil restrictions of scalars.} Let $i_1:R_1\hookrightarrow R$ be a finite, flat $\dbZ$-monomorphism. Let $\Lie_{R_1}(S)$ be $\Lie(S)$ but viewed as an $R_1$-Lie algebra. If $R_1=\dbF_p$, then we often view $\Lie_{\dbF_p}(S)$ just as an abelian group; if also $R=\dbF_p$, then we often drop the lower right index $\dbF_p$. Let $\Res_{R/R_1} S$ be the affine group scheme over $R_1$ obtained from $S$ through the Weil restriction of scalars (see [2, 7.6] and [7, 1.5]). So $\Res_{R/R_1} S$ is defined by the functorial group identification 
$$\Hom_{\Spec(R_1)}(Y,\Res_{R/R_1} S)=\Hom_{\Spec(R)}(Y\times_{\Spec(R_1)} \Spec(R),S),\leqno (1)$$
where $Y$ is an arbitrary $\Spec(R_1)$-scheme. Based on 2.2 and (1) we get a canonical and functorial identification 
$$\Lie(\Res_{R/R_1} S)=\Lie_{R_1}(S).\leqno (2)$$
If $i_2:R_2\hookrightarrow R_1$ is a second finite, flat $\dbZ$-monomorphism, then we have a canonical and functorial identification 
$$\Res_{R_1/R_2}\Res_{R/R_1} S=\Res_{R/R_2} S.$$ 
\noindent
{\bf 2.3.1. Proposition.} {\it We assume there is a finite subgroup $C$ of $\text{Aut}_{R_1}(R)$ such that we have an $R$-isomorphism $R\otimes_{R_1} R=\prod_{c\in C} R\otimes_R {}_c R$. Then we have
$$(\Res_{R/R_1} S)_R=\prod_{c\in C} S\times_{\Spec(R)} {}_c\Spec(R).\leqno (3)$$
\indent
If $S_1$ is a closed (resp. smooth, flat, semisimple or reductive) subgroup of $S$, then $\Res_{R/R_1} S_1$ is a closed (resp. smooth, flat, semisimple or reductive) subgroup of $\Res_{R/R_1} S$. If $S$ is a reductive group scheme, then any maximal torus (resp. Borel subgroup) $J$ of $\Res_{R/R_1} S$ is of the form $\Res_{R/R_1} T$ (resp. $\Res_{R/R_1} B$), where $T$ (resp. $B$) is a maximal torus (resp. Borel subgroup) of $S$.}

\medskip
\proof
Formula (3) follows from (1) and our hypothesis. See [2, Prop. 5 of 7.6] for the closed subgroup, smooth or flat part. Let now $S$ be a reductive group scheme. We know that $\Res_{R/R_1} S$ is affine, flat and smooth. So in order to show that it is a reductive group scheme, it is enough to show that its geometric fibres are so. So it is enough to show that the fibres of $(\Res_{R/R_1} S)_R$ are reductive group schemes. But this is so as $(\Res_{R/R_1} S)_R$ is a reductive group scheme, cf. (3). Similarly we argue that $\Res_{R/R_1} S$ is semisimple if $S$ is so. 

To check the last part, we first remark that $J_R$ is a maximal torus (resp. Borel subgroup) of $(\Res_{R/R_1} S)_R$. So based on (3), it is uniquely determined by the projection $T$ (resp. $B$) of $J_R$ on the factor $S=S\times_{\Spec(R)} {}_{1_R}\Spec(R)$ of $(\Res_{R/R_1} S)_R$. So $J=\Res_{R/R_1} T$ (resp. $J=\Res_{R/R_1} B$). This ends the proof. 

\medskip
The next structure Theorem (see [29, 3.1.2]) will play important roles in \S3 and \S4.

\medskip\noindent  
{\bf 2.3.2. Theorem.} {\it We assume $R$ is a field. Then any adjoint (resp. s.c.) group $F$ over $R$ is isomorphic to a product of adjoint (resp. s.c.) groups of the form $\Res_{\tilde R/R} \tilde F$, where $\tilde R$ is a separable finite field extension of $R$ and $\tilde F$ is an adjoint (resp. s.c.) group over $\tilde R$ having an a.s. adjoint. So if $F^{\ad}$ is simple, then $F$ is of isotypic Dynkin type.}

                                                                               \bigskip
\noindent                                             
{\boldsectionfont \S 3. Adjoint representations over finite fields}

\bigskip
Let $H$ be a semisimple group over $k=\dbF_{q}$. In \S 3 we study $H$, $\Lie(H)$ and $AD_H$. In 3.3 we recall all cases when $\Lie(H)$ is simple, i.e. it has no proper ideal defined over $k$. These cases are precisely the ones when $AD_H$ is irreducible, cf. 3.4. In 3.5 to 3.11 we work in the context of Weil restrictions. In 3.1, 3.2 and 3.5 to 3.8 we include preliminary material on $H$ and $\Lie(H)$. See 3.7.1 for a refinement of 3.4 obtained by working with 
$$\grL_{H}:=\im(\Lie(H^{\sc})\to\Lie(H))$$ 
instead of $\Lie(H)$. In 3.9 we list semisimple subgroups of $H$ to be used in 3.10, 3.11 and \S 4. In 3.10 and 3.11 we study subrepresentations of $AD_H$; these sections extend [13]. 

\medskip\noindent
{\bf 3.1. On $Z(H)$.} The finite group scheme $Z(H)$ is of multiplicative type, cf. [27, Vol. III, 4.1.7 of p. 173]. So the finite group $Z(H)(k)$ is of order prime to $p$.  If $p|o(H)$, then $\Lie(Z(H))$ is a proper ideal of $\Lie(H)$ normalized by $H$ and so $AD_H$ is reducible. If $g.c.d.(p,o(H))=1$, then the central isogeny $H\to H^{\ad}$ is \'etale and so by identifying the tangent spaces in the origin we get $\Lie(H^{\ad})=\Lie(H)$. Similarly, if  $g.c.d.(p,c(H))=1$ we have $\Lie(H^{\sc})=\Lie(H)$.

\medskip\noindent
{\bf 3.2. Lemma.} {\it We assume that $H^{\ad}$ is a.s. and $\Lie(H^{\ad})$ is not a.s. Then $\Lie(H^{\ad})$ has a proper, characteristic ideal and so $AD_H$ is reducible.}

\medskip
\proof
If $p\Ge 5$, or if $p=3$ and $H^{\ad}$ is not of $G_2$ Lie type, or if $p=2$ and $H^{\ad}$ is not of $F_4$ Lie type, then $[\Lie(H^{\ad}),\Lie(H^{\ad})]$ is a characteristic ideal of $\Lie(H^{\ad})$ of codimension 1 or 2 (cf. [18, 0.13] applied to $H^{\ad}_{\kbar}$). In the excluded two cases the group $H^{\ad}$ is split and (cf. loc. cit.) $\Lie(H^{\ad})$ has a unique proper ideal: it corresponds to root vectors of short roots. So $\Lie(H^{\ad})$ has a proper, characteristic ideal and so $AD_H$ is reducible. 

\medskip\noindent
{\bf 3.3. Proposition.} {\it The Lie algebra $\Lie(H)$ is not simple iff one of the following four conditions holds:

\medskip
i) $H^{\ad}$ is not simple;

\smallskip
ii) $H^{\ad}$ is simple of isotypic $A_{pn-1}$ Dynkin type;

\smallskip
iii) $p=2$ and $H^{\ad}$ is simple of an isotypic Dynkin type belonging to the union $\{B_n,C_n|n\in\dbN\}\cup\{D_n|n\Ge 3\}\cup\{E_7,F_4\}$;

\smallskip
iv)  $p=3$ and $H^{\ad}$ is simple of isotypic $E_6$ or $G_2$ Dynkin type.}

\medskip
\proof
If i) to iv) do not hold, then $g.c.d.(p,o(H))=1$ and so $\Lie(H)=\Lie(H^{\ad})$. If $\Lie(H)$ is simple, then $g.c.d.(p,o(H))=1$ and so $\Lie(H)=\Lie(H^{\ad})$ (cf. 3.1) and $H$ is simple (cf. 2.3.2). So it suffices to prove the Proposition under the extra assumption that $H$ is adjoint and simple. Let $k_1$ be such that $H=\Res_{k_1/k} H_1$, where $H_1$ is an a.s. adjoint group over $k_1$ (cf. 2.3.2). We have $\Lie(H)=\Lie_k(H_1)$, cf. (2). We show that $\Lie(H_1)$ is simple iff $\Lie(H)$ is simple. The ``if" part is obvious. To check the ``only if" part let $I$ be a non-zero ideal of $\Lie(H)$. So $I\otimes_k k_1$ is a non-zero ideal of 
$$\Lie(H)\otimes_k k_1=\Lie_k(H_1)\otimes_k k_1\tilde\to\oplus_{\tau\in\Gal(k_1/k)} \Lie(H_1)\otimes_{k_1} {}_{\tau} k_1\leqno (4)$$ 
and so a direct sum of some of these direct summands indexed by elements of $\Gal(k_1/k)$. But the only non-zero such direct sum defined over $k$ is $\Lie_k(H_1)\otimes_k k_1$. So $I=\Lie(H)$. 

Also, $\Lie(H_1)$ is a.s. iff it is simple (cf. 3.2 for the if part). So $\Lie(H)$ is simple iff $\Lie(H_1)$ is a.s. So the Proposition follows from [18, 0.13].

\medskip\noindent
{\bf 3.4. Proposition.} {\it Let $H$ be a semisimple group over $k$. Then $AD_H$ is irreducible iff $\Lie(H)$ is simple (i.e. iff none of the four conditions of 3.3 holds).}

\medskip
\proof 
We assume $AD_H$ is irreducible. So $g.c.d.(p,o(H))=1$ and $\Lie(H)=\Lie(H^{\ad})$, cf. 3.1. Moreover $H^{\ad}$ is simple, cf. 2.3.2. So it suffices to show that $\Lie(H)$ is simple under the extra assumption that $H^{\ad}$ is a.s., cf. the proof of 3.3. So $\Lie(H)$ is a.s., cf. 3.2. 

We now assume $\Lie(H)$ is simple. So $H^{\ad}$ is simple and $g.c.d.(p,o(H^{\sc}))=1$, cf. 3.3. So to show that $AD_H$ is irreducible we can assume that $H=\Res_{k_1/k} H_1$, with $H_1$ a s.c. semisimple group over $k_1$ having an a.s. adjoint. So $\Lie(H_1)=\Lie(H_1^{\ad})$ are a.s., cf. 3.3. Let $DT$ be the Dynkin type of $H_1$. The maximal weight $\varpi$ of the adjoint representation of the complex, simple Lie algebra of $DT$ Dynkin type is the maximal root. So if $DT$ is $A_n$ (resp. is $B_n$ with $n\Ge 3$, $C_n$ with $n\Ge 2$, $D_n$ with $n\Ge 4$, $E_6$, $E_7$, $E_8$, $F_4$ or $G_2$), then with the notations of [5, planche I to IX] $\varpi$ is $\varpi_1+\varpi_n$ (resp. is $\varpi_2$, $2\varpi_1$, $\varpi_2$, $\varpi_2$, $\varpi_1$, $\varpi_8$, $\varpi_1$ or $\varpi_2$). From Curtis and Steinberg theory and the last two sentences we get that $AD_{H_1}$ is absolutely irreducible. Let $I$ be a non-trivial, irreducible subrepresentation of $AD_H$. The extension of $AD_H$ to $k_1$ is a direct sum of $[k_1:k]$ copies of $AD_{H_1}$, cf. (4). So due to the absolute irreducibility of $AD_{H_1}$, $I\otimes_k k_1$ is a $\Lie(H_1)$-submodule of $\Lie(H)\otimes_k k_1$. So as $(I\otimes_k k_1)\cap \Lie(H)=I$, $I$ is an ideal of $\Lie(H)$. So $I=\Lie(H)$. This ends the proof.

\medskip\noindent
{\bf 3.5. The basic setting.} Until \S 4 we assume that 
$$H=\Res_{k_1/k}H_1,$$ 
where $H_1$ is a semisimple group over $k_1$ such that $H_1^{\ad}$ is a.s. So $\Lie(H)=\Lie_k(H_1)$. Let 
$$DT$$ 
be the Dynkin type of $H_1$. Let $k_2$ and $k_3$ be the quadratic and respectively the cubic field extension of $k_1$. The group $H_1$ has a Borel subgroup $B_1$ (this is part of Lang theorem). Let $T_1$ be a maximal torus of $B_1$, cf. [4, (i) of Theorem 18.2]. The group $H_1$ is split iff $T_1$ is. So there is a smallest field extension $l$ of $k_1$ such that $H_{1l}$ is split. Let
$$\Lie(H_{1l})=\Lie(T_{1l})\bigoplus\oplus_{\alpha\in\Phi} {\grg}_{\alpha}\leqno (5)$$
be the Weyl decomposition of $\Lie(H_{1l})$ with respect to $T_{1l}$. So each ${\grg}_{\alpha}$ is a 1 dimensional $l$-vector space normalized by $T_{1l}$ and $\Phi$ is an irreducible root system of characters of $T_{1l}$. Let $\Phi^+:=\{\alpha\in\Phi|\grg_{\alpha}\subset\Lie(B_{1l})\}$. Let $\Delta$ be the basis of $\Phi$ included in $\Phi^+$. We use [5, planches I to IX] to denote the elements of $\Delta$, $\Phi^+$ and $\Phi$; so $\Delta=\{\alpha_1,...,\alpha_{|\Delta|}\}$, where $|\Delta|$ is the number of elements of $DT$. Warning: if $l=k_1$ (resp. $l=k$), then in connection to (5) we often drop the lower right index $l$ (resp. $1l$ and so also $T_1$ becomes $T$).

We denote also by $\Delta$ the Dynkin diagram of $DT$ defined by the elements of $\Delta$. Let $\phi$ be the Frobenius automorphism of $\bar l$ fixing $k_1$. It acts on $\Phi$ via the rule: if $\alpha\in\Phi$, then we have $\grg_{\phi(\alpha)}=(1_{\Lie(H_1)}\otimes\phi)(\grg_{\alpha})$. Let $e\in\dbN$ be the order of the permutation of $\Delta$ defined by $\phi$. For $m\in\dbN$, $\phi^m$ fixes $\Lie(T_{1l})$ iff $e|m$. So $[l:k_1]=e$. But $e$ is the order of an element of $\text{Aut}(\Delta)$ and so $e\in\{1,2,3\}$. More precisely, we have  $l=k_1$ if $\in\{B_n,C_n|n\Ge 1\}\cup\{E_7,E_8,F_4,G_2\}$, $l\in\{k_1,k_2\}$ if $DT\in\{A_n|n\Ge 2\}\cup\{D_n|n\Ge 5\}\cup\{E_6\}$ and $l\in\{k_1,k_2,k_3\}$ if $LT=D_4$. If $e=1$ let $LT:=DT$. If $e\in\{2,3\}$ let $LT:={}^eDT$. We recall that the pair $(\Delta,\phi)$ or just 
$$LT$$ 
is called the Lie type of $H_1$ and that $H_1^{\ad}$ and $H_1^{\sc}$ are uniquely determined by it (see [29]).

For $\alpha\in\Phi$ let $H(|\alpha|)$ be the semisimple subgroup of $H_{1l}$ generated by the two $\dbG_a$ subgroups $\dbG_{a,\alpha}$ and $\dbG_{a,-\alpha}$ of $H_{1l}$ normalized by $T_{1l}$ and having ${\grg}_{\alpha}$ and respectively ${\grg}_{-\alpha}$ as their Lie algebras. So $H(|\alpha|)^{\ad}$ is a $PGL_2$ group and $H_{1l}$ is generated by all these $H(|\alpha|)$'s. All these can be deduced from [4, \S 13 of Ch. IV] via descent from $\lbar$ to $l$.

\medskip\noindent
{\bf 3.6. Proposition.} {\it Let $\beta\in\Phi\setminus\{-\alpha,\alpha\}$. Let $\Psi_1(\alpha,\beta):=\{i\alpha+\beta|i\in\dbN,\, i\alpha+\beta\in\Phi\}$. If $\alpha+\beta\notin\Phi$ let $\grg_{\alpha+\beta}:=\{0\}$. Let $g\in\dbG_{a,\alpha}(l)$ be a non-identity element. Let $x\in\grg_{\beta}\setminus\{0\}$. 

\medskip
{\bf a)} We have $AD_{H_{1l}}(g)(x)-x\in\oplus_{\delta\in\Psi_1(\alpha,\beta)} \grg_{\delta}$.

{\bf b)} We always have an inclusion $[\grg_{\alpha},\grg_{\beta}]\subset\grg_{\alpha+\beta}$. It is not an equality iff one of the following disjoint three conditions holds:

\medskip
i) $p=3$, $DT=G_2$, and $\alpha$ and $\beta$ are short roots forming an angle of $60^o$;

\smallskip
ii) $p=2$, $DT=G_2$, and $\alpha$ and $\beta$ are short roots forming an angle of $120^o$;

\smallskip 
iii) $p=2$, $DT\in\{B_n,C_n|n\Ge 2\}\cup\{F_4\}$, $\alpha$ and $\beta$ are short, $\alpha\perp\beta$ and $\alpha+\beta\in\Phi$.

\medskip
{\bf c)} If none of the three conditions of b) holds and if $\alpha+\beta\in\Phi$, then the component of $AD_{H_{1l}}(g)(x)-x$ in $\grg_{\alpha+\beta}$ is non-zero.}

\medskip
\proof
Let $\kappa(\alpha,\beta):=\Phi\cap\{i\alpha+j\beta|(i,j)\in\dbN^2\}$. It is known that there are isomorphisms $x_{\delta}:\dbG_a\arrowsim\dbG_{a,\delta}$ ($\delta\in\kappa(\alpha,\beta)$), an ordering on $\dbN^2$, and integers $m_{\alpha,\beta;i,j}$ ($(i,j)\in\dbN^2$) such that for all $s$, $t\in\dbG_a(\bar l)=\bar l$ we have (cf. [27, Vol. III, p. 321] or [8, 4.2])
$$x_{\alpha}(s)x_{\beta}(t)x_{\alpha}(s^{-1})x_{\beta}(t^{-1})=\prod_{(i,j)\in\dbN^2,\,i\alpha+j\beta\in\kappa(\alpha,\beta)} \,x_{i\alpha+j\beta}(m_{\alpha,\beta;i,j}s^it^j).\leqno (6)$$
 Taking the derivative in (6) with respect to $t$ at the identity element, we get that the difference $AD_{H_{1\bar l}}(x_{\alpha}(s))(x)-x$ belongs to $(\oplus_{i\in\dbN,\,i\alpha+\beta\in\Psi_1(\alpha,\beta)} \grg_{i\alpha+\beta})\otimes_l\otimes \bar l$ and its component in $\grg_{i\alpha+\beta}\otimes_l\otimes \bar l$ is a fixed generator of $\grg_{i\alpha+\beta}\otimes_l\otimes \bar l$ times $m_{\alpha,\beta;i,1}s^i$ (cf. [4, 3.16 of Ch. 1]). Similarly, taking the derivative in (6) with respect to $s$ and $t$ at the identity element, we get that $[\grg_{\alpha},\grg_{\beta}]=m_{\alpha,\beta;1,1}\grg_{\alpha+\beta}\subset \grg_{\alpha+\beta}$. But $m_{\alpha,\beta;1,1}=\pm(m+1)$, where $m\in\dbZ$ is the greatest integer such that $\beta-m\alpha\in\Phi$ (see [27, Vol. III, p. 321] or [8, 4.2 and 4.3]). So the second part of b) and c) follow from the fact that we have $p|m+1$ iff one of the three conditions of b) holds (see loc. cit.). This ends the proof.

\medskip
Let $B_1^{\text{opp}}$ be the Borel subgroup of $H_1$ which is the opposite of $B_1$ with respect to $T_1$. Let $U_1$ and $U_1^{\text{opp}}$ be the unipotent radicals of $B_1$ and respectively of $B_1^{\text{opp}}$. Let $L_H$ be the Lie subalgebra of $\Lie(H)=\Lie_k(H_1)$ generated by $\Lie_k(U_1)$ and $\Lie_k(U_1^{\text{opp}})$. 

\medskip\noindent
{\bf 3.7. Proposition.} {\it {\bf 1)} The Lie algebra $L_H$ is an ideal of $\Lie(H)$ containing $[\Lie(H),\Lie(H)]$. If $p>2$ and $l=k_1$, then $L_H=\sum_{\alpha\in\Phi^+} \Lie_k(H(|\alpha|))$. 

\smallskip
{\bf 2)} Always $L_H=\grL_H$. Also, if $p>2$ or if $p=2$ and $DT\notin\{C_n|n\in\dbN\}$, then $L_H=[\Lie(H),\Lie(H)]$. 

\smallskip
{\bf 3)} We have  $\grL_{H^{\ad}}=L_{H^{\ad}}=[\Lie(H^{\ad}),\Lie(H^{\ad})]$.}

\medskip
\proof
Tensoring with $l$ we can assume $k=k_1=l$; so $H=H_1$ is split. The first part of 1) is obvious. If $p\Ge 3$, then $\Lie(H(|\alpha|))=\grg_{\alpha}\oplus\grg_{-\alpha}\oplus [\grg_{\alpha},\grg_{-\alpha}]$. So the last part of 1) follows from the inclusions of 3.6 b). We now prove 2). Both $U_1$ and $U_1^{\text{opp}}$ are naturally identified with unipotent subgroups of $H_1^{\sc}$ and so $L_{\tilde H}$ makes sense if $\tilde H\to H$ is an arbitrary central isogeny. Always $L_{\tilde H}$ surjects onto $L_H$ and so to prove 2) we can assume that $H$ is s.c. If $p>2$ or if $p=2$ and $DT\notin\{C_n|n\in\dbN\}$, then $\Lie(H)=[\Lie(H),\Lie(H)]$ (cf. [18, 0.13]) and so 2) follows from 1). So we are left to show that $L_H=\Lie(H)$ if $p=2$ and $DT=C_n$. But then $H$ has a subgroup $SL(H)$ which is a product of $SL_2$ groups and has $T$ as a maximal torus; so $\Lie(T)\subset L_H$ and so $L_H=\Lie(H)$. This proves 2). 
 
Based on 2) and its proof, its suffices to prove 3) under the extra assumptions that $p=2$ and $DT=C_n$. As $L_{H^{\ad}}=\grL_{H^{\ad}}$, the ideal $L_{H^{\ad}}$ of $\Lie(H^{\ad})$ is of codimension 1 and contains $[\Lie(H^{\ad}),\Lie(H^{\ad})]$. By reasons of dimensions we get $L_{H^{\ad}}=[\Lie(H^{\ad}),\Lie(H^{\ad})]$ (see [18, 0.13]). This ends the proof.

\medskip\noindent
{\bf 3.7.1. Corollary.} {\it We assume $g.c.d.(p,o(H))=1$. If $p=3$  we also assume that $DT\neq G_2$ and if $p=2$ we also assume that $DT\notin\{B_n,C_n|n\Ge 2\}\cup\{F_4\}$. Then the natural representation of $H^{\sc}(k)=H_1^{\sc}(k_1)$ on $\grL_{H^{\ad}}=\grL_H=[\Lie(H),\Lie(H)]$ is irreducible.}

\medskip
\proof 
The  $H_1$-module $\grL_{H_1}$ is a.s. Argument: the case $DT=A_1$ is trivial and the case $DT\neq A_1$ is a consequence of the fact that $\grL_{H_1}$ is an a.s. $\Lie(H_1)$-module (see [18, 0.13]). The maximal weight of the representation of $H_{1l}$ on $\grL_{H_{1l}}$ is $\varpi$ of the proof of 3.4. The rest of the proof is entirely the same as the proof of 3.4.

\medskip
We come back to (5). Let $\alpha\in\Phi$. Let $T(|\alpha|)$ be the maximal subtorus of $T_{1l}$ centralizing $\grg_{\alpha}$. So $[\grg_{\alpha},\Lie(T(|\alpha|))]=\{0\}$. We have a short exact sequence
$$0\to T(|\alpha|)\to TH(|\alpha|)\to\tilde H(|\alpha|)\to 0,\leqno (7)$$
where $TH(|\alpha|)$ is the reductive subgroup of $H_{1l}$ generated by $T(|\alpha|)$ and $H(|\alpha|)$ and where $\tilde H(|\alpha|)$ is either $H(|\alpha|)$ or its adjoint.

\medskip\noindent
{\bf 3.8. Lemma.} {\it The group $H(|\alpha|)$ is a $PGL_2$ group iff $DT=B_n$, $H_1$ is adjoint and $\alpha$ is short. Also, if $H_1$ is adjoint, then $\tilde H(|\alpha|)=H(|\alpha|)^{\ad}$. 
}

\medskip
\proof  
For $\beta\in\Phi$ let $ST(\alpha,\beta):=\{i\alpha+\beta|i\in\dbZ\}\cap\Phi$ be the $\alpha$-string through $\beta$. It is of the form $\{i\alpha+\beta|i\in\{-s,-s+1,...,t\}\}$, for some $s$, $t\in\dbN\cup\{0\}$ with $s+t\Le 3$ (see [15, p. 45]). The set $\Phi\setminus\{\alpha,-\alpha\}$ is a disjoint union of $\alpha$-strings. So let $\Phi(\alpha)$ be a subset of $\Phi\setminus\{\alpha,-\alpha\}$ such that we have a disjoint union
$\Phi\setminus\{\alpha,-\alpha\}=\cup_{\beta\in \Phi(\alpha)} ST(\alpha,\beta)$. For $\beta\in\Phi(\alpha)$ let $V_{\alpha,\beta}:=\oplus_{\delta\in ST(\alpha,\beta)} \grg_{\delta}$. It is an $H(|\alpha|)$-module, cf. 3.6 a). We get a direct sum decomposition of $H(|\alpha|)$-modules $\Lie(H_{1l})=\Lie(TH(|\alpha|))\bigoplus\oplus_{\beta\in\Phi(\alpha)} V_{\alpha,\beta}$.

If $H(|\alpha|)$ is a $PGL_2$ group, then $\dim_{l}(V_{\alpha,\beta})\in\{1,3\}$, $\forall\beta\in\Phi(\alpha)$, and the converse holds if $H_1$ is adjoint. Using [5, planches I to IX] and [15, Lemma C of p. 53] we easily get that $\dim_l(V_{\alpha,\beta})\in\{1,3\}$, $\forall\beta\in\Phi(\alpha)$, iff $DT=B_n$ and $\alpha$ is short. But if $H_1$ is s.c., then $H(|\alpha|)$ is an $SL_2$ group as we can check by reduction to the $B_1$ and $B_2=C_2$ Dynkin types. This proves the first part. 

If $H_{1}$ is adjoint, then the subgroup of $T_{1l}$ fixing $\grg_{\alpha}$ is a subtorus of $T_{1l}$ of codimension 1 and so it is $T(|\alpha|)$. So $Z(H(|\alpha|))$ is a subgroup of $T(|\alpha|)$ and so $\tilde H(|\alpha|)=H(|\alpha|)^{\ad}$. This ends the proof.

\medskip\noindent
{\bf 3.8.1. Remark.} If $(\Psi_1(\alpha,\beta)\cup \Psi_1(-\alpha,\beta))\neq\emptyset$, then $\Phi\cap (\dbZ\alpha\oplus\dbZ\beta)$ is the root system of a Dynkin type $DT_2:=DT(\alpha,\beta)\in\{A_2,B_2,G_2\}$. If $DT_2=G_2$, then $DT=G_2$. If $DT_2=B_2$, then $DT\in\{B_n,C_n|n\Ge 2\}\cup\{F_4\}$. Let $x\in\grg_{\alpha}$. If $DT_2=A_2$, then $\Psi_1(\alpha,\beta)$ has at most 1 element and so $\ad(x)^2$ annihilates $V_{\alpha,\beta}$ (cf. 3.6 a)). Similarly if $DT_2=B_2$ (resp. $DT_2=G_2$), then $\Psi_1(\alpha,\beta)$ has at most 2 (resp. 3) elements and so $\ad(x)^3$ (resp. $\ad(x)^4$) annihilates $V_{\alpha,\beta}$.

\medskip\noindent
{\bf 3.9. Subgroups.} We list semisimple subgroups of $H_1$ normalized by $T_1$. By taking the $\Res_{k_1/k}$ of them we get semisimple subgroups of $H$. So except for the first two Cases below (which are samples) we will just mention semisimple subgroups of $H_1$ and not of $H$. 

Let $\Phi_0$ be a subset of $\Phi$ having the following three properties:

\medskip
(1) it is invariant under $\phi$ of 3.5;

(2) it is stable under additions, i.e. we have $\Phi\cap\{\alpha+\beta|\alpha,\beta\in\Phi_0\}\subset\Phi_0$;

(3) it is symmetric, i.e. we have $\Phi_0=-\Phi_0$.

\medskip 
Let $H_{0l}$ be the subgroup of $H_{1l}$ generated by $H(|\alpha|)$, where $\alpha\in\Phi_0$. It is invariant under $\phi$ and $T_{1l}$ and so it the extension to $l$ of a subgroup $H_0$ of $H_1$ normalized by $T_1$. The subgroup $TH_0$ of $H_1$ generated by $T_1$ and $H_0$ is reductive, cf. [27, Vol. III, 5.4.7 and 5.10.1 of Exp. XXII] applied to $TH_{0l}$. But $H_0$ is a subgroup of $(TH_0)^{\der}$ and contains the unipotent radical of any Borel subgroup of $TH_0$ normalized by $T_1$. So $H_0=(TH_0)^{\der}$ and so $H_0$ is semisimple. A subset $\Phi_0^+$ of $\Phi_0\cap\Phi^+$ is said to generate $\Phi_0$, if each element of $\Phi_0\cap\Phi^+$ is a linear combination with non-negative integer coefficients of elements of $\Phi_0^+$. We refer to $H_0$ as the semisimple subgroup of $H_1$ associated to $\Phi_0$ or generated by $\Phi_0^+$. 

In this paragraph we assume $l=k_1$. If $DT\notin\{A_1,B_2,G_2\}$, then we can choose $\Phi_0$ such that $H_0$ is of $A_2$ Lie type. If $DT\notin\{A_1,A_2,B_2,B_3,C_3,D_4,F_4,G_2\}$ (resp. $DT\in\{B_n,C_n|n\Ge 2\}\cup\{F_4\}$), then we can choose $\Phi_0$ such that $H_0$ is of $A_3$ (resp. $C_2$) Lie type. We assume now that $DT\in\{B_n,C_n|n\Ge 2\}\cup\{F_4,G_2\}$. Let $H_1^{\text{long}}$ be the semisimple subgroup of $H_1$ generated by $\{\alpha\in\Phi^+|\alpha\,\text{is long}\}$. Warning: the group $H_1^{\text{long}}$ as well as most of the below semisimple subgroups of $H_1$ are not generated by subsets of $\Delta$. If $n\Ge 3$ and $DT=B_n$, then $H_1^{\text{long}}$ is of $D_n$ Lie type. If $n\Ge 2$, $H_1$ is s.c. and $DT=C_n$, then  $H_1^{\text{long}}$ is an $SL_2^n$ group. If $DT=F_4$ (resp. $DT=G_2$), then $H_1^{\text{long}}$ is of $A_3$ (resp. $A_2$) Lie type. These Lie types can be read out from [5, planches II, III, VIII and IX].

Until 3.10 we use (5) with $l\in\{k_2,k_3\}$. If $\varpi$ is as in the proof of 3.4, then the subgroup $H(|\varpi|)$ of $H_{1l}$ is the pull back of an $SL_2$ subgroup of $H_1$ (cf. also 3.8 applied to $H_{1l}$). If $LT={}^2A_{n+4}$ and $\Phi_0^+=\{\alpha_1,\alpha_2,\alpha_{n+3},\alpha_{n+4}\}$, then $H_0$ is isogenous to the $\Res_{k_2/k_1}$ of an $SL_3$ group. If $LT={}^2A_4$ and $\Phi_0^+=\{\alpha_1,\alpha_4\}$, then $H_0$ is the $\Res_{k_2/k_1}$ of an $SL_2$ group. We assume now that $LT={}^2D_{n+3}$. If $n=1$ we assume $\phi$ permutes $\alpha_3$ and $\alpha_4$. If $n>1$ (resp. $n=1$) and $\Phi_0^+=\{\alpha_1,\alpha_2,\alpha_3\}$ (resp. $\Phi_0^+=\{\alpha_1,\alpha_2,\alpha_2+\alpha_3+\alpha_4\}$), then $H_0^{\ad}$ is a $PGL_4$ group. To introduce other semisimple subgroups of $H_1$ we consider five Cases. 

\medskip\noindent
{\bf Case 1: $LT={}^2A_{2n+1}$.} Let $H_1(n+1)$ be the $SL_2$ subgroup of $H_1$ such that $H_1(n+1)_l=H(|\alpha_{n+1}|)$. So $H(n+1):=\Res_{k_1/k} H_1(n+1)$ is a  subgroup of $H$. Let $i\in\{1,...,n\}$. Let $\alpha^1_i:=\sum_{j=i}^{n+1} \alpha_j$ and $\alpha^2_i:=\sum_{j=n+1}^{2n+2-i} \alpha_j$. The subgroups $H(|\alpha^1_i|)$ and $H(|\alpha^2_i|)$ of $H_{1k_2}$ are $SL_2$ groups and commute. We first assume that either $n>1$ or $H_1=H_1^{\sc}$. The product $H(|\alpha^1_i|)\times_{k_2} H(|\alpha^2_i|)$ is naturally a subgroup of $H_{1k_2}$, invariant under the action of $\Gal(k_2/k_1)$ on $H_{1k_2}$. So it is the extension to $k_2$ of a subgroup $H_1(i)$ of $H_1$. As $\Gal(k_2/k_1)$ takes $\Lie(H(|\alpha^1_i|))$ into $\Lie(H(|\alpha^2_i|))$, the group $H_1(i)$ is the $\Res_{k_2/k_1}$ of an $SL_2$ group. So $H(i):=\Res_{k_1/k} H_1(i)$ is a subgroup of $H$ which is the $\Res_{k_2/k}$ of an $SL_2$ group. If $n=1$ and $H_1\neq H_1^{\sc}$, then let $H_1(1)$ be the semisimple subgroup of $H_1$ such that we have a natural central isogeny $H(|\alpha^1_1|)\times_{k_2} H(|\alpha^2_1|)\to H_1(1)_{k_2}$ of degree 2; the adjoint of the subgroup $H(1):=\Res_{k_1/k} H_1(1)$ of $H$ is the $\Res_{k_2/k}$ of a $PGL_2$ group.

\medskip\noindent
{\bf Case 2: $LT={}^2A_{2n}$.} 
For $i\in\{1,...,n-1\}$ let $\alpha^1_i:=\sum_{j=i}^{n+1} \alpha_j$ and $\alpha^2_i:=\sum_{j=n}^{2n+1-i} \alpha_j$. Let $H_1(i)$ and $H_1(n+1)$ be subgroups of $H_1$ such that $H_1(i)_{k_2}=H(|\alpha^1_i|)\times_{k_2} H(|\alpha^2_i|)$ and $H_1(n+1)_{k_2}=H(|\alpha_n+\alpha_{n+1}|)$. The adjoint of the semisimple subgroup $H_1(n)$ of $H_1$ generated by $\{\alpha_n,\alpha_{n+1}\}$ is a $PGU_3$ group. So for $i\in\{1,...,n+1\}$, the group $H(i):=\Res_{k_1/k} H_1(i)$ is a  subgroup of $H$. If $i\in\{1,...,n-1\}$ (resp. if $i=n+1$), then $H(i)$ is the $\Res_{k_2/k}$ (resp. $\Res_{k_1/k}$) of an $SL_2$ group.

\medskip\noindent
{\bf Case 3: $LT={}^2E_6$.} Let $H_1(1)^\prime$ (resp. $H_1(4)$) be the  semisimple subgroup of $H_1$ generated by $\{\alpha_1,\alpha_3,\alpha_4,\alpha_5,\alpha_6\}$ (resp. $\{\alpha_2\}$). So $H_1(1)^\prime$ is of ${}^2A_5$ Lie type and $H_1(4)$ is an $SL_2$ group, cf. 3.8. We apply Case 1 to $H_1(1)^\prime$: we get semisimple subgroups $H_1(1)$, $H_1(2)$ and $H_1(3)$ of $H_1(1)^\prime$ and so also of $H_1$.

\medskip\noindent
{\bf Case 4: $LT={}^2D_{n+3}$.} If $n=1$ we assume $\phi$ permutes $\alpha_3$ and $\alpha_4$. For $i\in\{1,...,n+2\}$ let $H_1(i)$ be the semisimple subgroup of $H_1$ which over $k_2$ is generated by $H(|\alpha|)$, with $\alpha\in\{\sum_{j=i}^{n+2}\alpha_j,\alpha_{n+3}+\sum_{j=i}^{n+1}\alpha_j\}$. So $H_1(i)^{\sc}$ is the $\Res_{k_2/k_1}$ of an $SL_2$ group.

\medskip\noindent
{\bf Case 5: $LT={}^3D_4$.} Let $H_1(1)$ (resp. $H_1(2)$) be the semisimple subgroup of $H_1$ which over $k_3$ is generated by $H(|\alpha|)$, with $\alpha\in\{\alpha_1+\alpha_2+\alpha_3,\alpha_1+\alpha_2+\alpha_4,\alpha_2+\alpha_3+\alpha_4\}$ (resp. with $\alpha=\alpha_1+\alpha_2+\alpha_3+\alpha_4$). So $H_1(2)$ is an $SL_2$ group (cf. 3.8) and $H_1(1)^{\sc}$ is the $\Res_{k_3/k_1}$ of an $SL_2$ group.

\medskip\noindent
{\bf 3.9.1. Lemma.} {\it We refer to Case 4 (resp. Case 5). If $H=H^{\sc}$, then each $H_1(i)$ (resp. then $H_1(1)$) are (resp. is) s.c.}

\medskip
\proof 
To check this we can work with $H_l$. Even better, based on the existence of Chevalley group schemes over $\dbZ$ we can work over $\dbC$. Based on [6, \S13 of Ch. VIII] we need to show that any standard (resp. standard composite) monomorphism $SO_4\hookrightarrow SO_{2n+6}$ (resp. $SL_2\times SO_4\hookrightarrow SO_4\times SO_4\hookrightarrow SO_8$) over $\dbC$, lifts to a monomorphism $\text{Spin}_4\hookrightarrow \text{Spin}_{2n+6}$ (resp. $SL_2^3=SL_2\times \text{Spin}_4\hookrightarrow \text{Spin}_8$). To check this last statement we can work just in the context of $SO_{2n+6}$. The restriction of the spin representation of $\Lie(SO_{2m})$ to $\Lie(SO_{2m-2})$ via a standard monomorphism $SO_{2m-2}\hookrightarrow SO_{2m}$ is a disjoint union of two spin representations of $\Lie(SO_{2m-2})$, $\forall m\in\dbN$, $m\ge 4$ (see the description of spin representations in loc. cit.). From this the statement on lifts follows. This ends the proof.

\medskip\noindent
{\bf 3.9.2. Lemma.} {\it The Lie subalgebra $LS_H$ of $\Lie_k(H_1)=\Lie(H)$ generated by all intersections $\Lie_k(H_1(i)\cap U_1)$'s and $\Lie_k(H_1(i)\cap U_1^{\text{opp}})$'s is $\grL_H=L_H$.}

\medskip
\proof
Tensoring with $l$ over $k$, this follows easily from 3.6 b).

\medskip\noindent
{\bf 3.9.3. Simple properties.} Let $TH_1(i)$ be the reductive subgroup of $H_1$ generated by $T_1$ and $H_1(i)$. Let $T_1(i)$ be the maximal torus of $Z(TH_1(i))$. We have a short exact sequence 
$$0\to\Lie(T_1(i))\to\Lie(TH_1(i))\to\Lie(\tilde H_1(i))\to 0,\leqno (8)$$ 
where $\tilde H_1(i)$ is the quotient of $H_1(i)$ by $Z(H_1(i))\cap T_1(i)$. If $H$ is adjoint, then as in the proof of 3.8 we argue that in fact $\tilde H_1(i)=H_1(i)^{\ad}$. Let $S_1(i):=\Lie(TH_1(i))$. Let
$$S_2(i):=\Lie(H_1)\cap (\oplus_{\alpha\in\Phi,\,\grg_{\alpha}\not\subset\Lie({H_1(i)}_l)} \grg_{\alpha}).$$
We have a direct sum decomposition $\Lie(H_1)=S_1(i)\oplus S_2(i)$ of $H_1(i)$-modules.
 
We assume now that $H_1(i)^{\ad}$ is the $\Res_{l/k_1}$ of a $PGL_2$ group. Let $x\in \Lie(U_1)\cap\Lie(H_1(i))$. In the first four Cases (resp. in Case 5) we can write $x=x_1+x_2$ (resp. $x=x_1+x_2+x_3$), where $x_i\in \grg_{\alpha_i}$ for some $\alpha_i\in\Phi^+$. This writing is unique up to ordering. The restriction of $\ad(x_i)^2$ to $S_2(i)$ is $0$, cf. the $DT_2=A_2$ part of 3.8.1. As $H_1(i)^{\ad}$ is the $\Res_{l/k_1}$ of a $PGL_2$ group we have $\ad(x_1)\ad(x_2)=\ad(x_2)\ad(x_1)$ (resp. $\ad(x_i)\ad(x_j)=\ad(x_j)\ad(x_i)$, $\forall i$, $j\in\{1,2,3\}$). Also in Case 5 we have $\ad(x_1)\ad(x_2)\ad(x_3)=0$. So 
$$h_x:=\sum_{s=0}^\infty{{\ad(x)^s}\over {s!}}=\sum_{s=0}^2 {{\ad(x)^s}\over {s!}}$$ 
is a well defined automorphism of $S_2(i)$. As in characteristic 0 we can identify $h_x$ with the automorphism of $S_2(i)$ defined by an element of $H_1(i)(k_1)$. 

We assume now that $LT={}^2A_{2n}$. So $H_1(n)$ is of ${}^2A_2$ Lie type. Let $x\in \Lie(H_1(n))\cap {\grg}_{\alpha_n}\oplus {\grg}_{\alpha_{n+1}}$. We write $x=x_1+x_2$, with $x_1\in {\grg}_{\alpha_n}$ and $x_2\in {\grg}_{\alpha_{n+1}}$. The restrictions of  $\ad(x_1)$ and $\ad(x_2)$ to $S_2(n)$ do not commute but $\ad(x_{i_1})\ad(x_{i_2})\ad(x_{i_3})$ restricted to $S_2(n)$ is $0$, $\forall i_1$, $i_2$, $i_3\in \{1,2\}$. So for $p\ge 3$ we can define $h_x$ as above.

\medskip\noindent
{\bf 3.10. Theorem.} {\it {\bf 1)} If $g.c.d.(p,o(H))=1$, then no 1 dimensional $k$-vector subspace of $\Lie(H)$ is normalized by $AD_H(\im(H^{\sc}(k)\to H(k)))$.}

\smallskip
{\it {\bf 2)} We assume that $p|o(H)$ and $g.c.d.(p,c(H))=1$. If $q=2$, we also assume that $H$ is not an $SL_2$ group. Then there is no proper $k$-vector subspace $V$ of $\Lie(H)$ normalized by $AD_H(\im(H^{\sc}(k)\to H(k)))$ and such that $\Lie(H)=V+\Lie(Z(H))$. 

\smallskip
{\bf 3)} We assume that $DT=D_{2n+2}$, that $o(H_1)=2$ and that $p=2$. Then there is no proper $k$-vector subspace $V_0$ of $[\Lie(H),\Lie(H)]$ normalized by $AD_H(\im(H^{\sc}(k)\to H(k)))$ and such that $[\Lie(H),\Lie(H)]=V_0+\Lie(Z(H))$.

\smallskip
{\bf 4)} We assume that $p$ divides both $c(H)$ and $o(H)$. If $p=2$ we also assume that $DT\neq D_{2n+2}$. Then the adjoint representation of $H^{\sc}(k)$ on $\Lie(H)$ is the direct sum of an irreducible, non-trivial representation and of a trivial representation.

\smallskip
{\bf 5)} If $p=2$ and $DT=C_n$ we assume that $H_1$ is adjoint. Then the $H_1(k_1)$-module $\Lie(H_1)/[\Lie(H_1),\Lie(H_1)]$ is trivial and of dimension over $k_1$ at most 2.}

\medskip
\proof
If $l\neq k_1$ we use the notations of 3.9. To prove 1) we can assume $H$ is adjoint and $k_1=k$. Let $u\in\Lie(H)$ be such that $ku$ is normalized by $AD_H(\im(H^{\sc}(k)\to H(k)))$. The case $DT=A_1$ is well known (we can assume that $q=p$ and cf. 3.4 that $p=2$). So we can assume $DT\neq A_1$. We first assume $l\neq k$. Let $x\in\Lie(H_1(i)\cap U_1)$, with $H_1(i)$ a semisimple subgroup of $H=H_1$ listed in the five Cases of 3.9. We write $u=u_1+u_2$, where $u_s\in S_s(i)$ ($s\in\{1,2\}$).
The group $H_1(i)(k_1)$ normalizes $ku_1$ and $ku_2$, cf. 3.9.3. 

We assume now that $DT\neq A_{2n}$. So $H_1(i)$ is of isotypic $A_1$ Dynkin type. We show that $[x,u_2]=0$. Let $a\in l\setminus\{0,1\}$. If $H_1(i)^{\ad}$ is a $PGL_2$ group, then the restriction of $\ad(x)^2$ to $S_2(i)$ is $0$  (cf. 3.8.1 applied to $H_{1l}$) and so the automorphism of $S_2(i)$ induced by $\sum_{s=0}^1\ad(x)^s$ is defined by an element of $H_1(i)(k)$. So $u_2+[x,u_2]\in ku_2$ and so as $\ad(x)$ is nilpotent we get $[x,u_2]=0$. If $H_1(i)^{\ad}$ is the $\Res_{l/k}$ of a $PGL_2$ group, then as $h_x$ and $h_{ax}$ (of 3.9.3) normalize $ku_2$ we similarly get that $[x,u_2]=0$. We now show that $[x,u_1]=0$. We use (8). If $p>2$, then from 3.4 we get that the image of $u_1$ in $\Lie(\tilde H_1(i))$ is $0$; so $[x,u_1]=0$. If $p=2$, then $\tilde H_1(i)(k_1)=H_1(i)^{\ad}(k_1)$ and as $\Lie(H_1(i)^{\ad})$ has no non-zero elements normalized by $H_1(i)^{\ad}(k_1)$ (cf. the case $DT=A_1$) we get $[x,u_1]=0$. 

So $[x,u]=0$. Similarly we get $[y,u]=0$, where $y\in\Lie(H_1(i)\cap U_1^{\text{opp}})$. So $u$ is annihilated by $\grL_H$, cf. 3.9.2. So $u$ as an element of $\Lie(H_l)$ annihilates ${\grg}_{\alpha}$, $\forall\alpha\in\Phi$. As $l\neq k$ none of the three conditions of 3.6 b) holds and so from 3.6 b) we get that the component of $u$ in ${\grg}_{\alpha}$ with respect to (5) is $0$, $\forall\alpha\in\Phi$. So $u\in\Lie(T_l)$. But as $H$ is adjoint, no element of $\Lie(T_l)$ annihilates ${\grg}_{\alpha}$, $\forall\alpha\in\Phi$. So $u=0$. 

We assume now that $DT=A_{2n}$. If $g.c.d.(p,2n+1)=1$, then the equality $u=0$ is implied by 3.4. So we can assume $p\Ge 3$. If $i\neq n$, then as above we get $[x,u]=0$. If $i=n$ we assume that $x$ is as in the end of 3.9.3. As $h_x$ and $h_{2x}$ normalize $ku_2$, as above we get $[x,u_2]=0$. So as in the previous paragraph to show that $u=0$ it is enough to show that $[x,u_1]=0$. To show this we can assume $n=1$ (so $\Lie(H)=S_1(1)$ and $u=u_1$). As $u$ is annihilated by $\pm\grg_{\alpha_1+\alpha_2}$ (cf. the case $i=n+1$) from 3.6 b) we get that the component of $u$ in $\grg_{\alpha}$ with respect to (5) is 0, $\forall\alpha\in\Phi$. So $u\in\Lie(T)$. Let $g\in U_1(k_1)\setminus \dbG_{a,\alpha_1+\alpha_2}(k_1)$. As $AD_H(g)(ku)\subset ku$, we get that $u$ is annihilated by ${\grg}_{\alpha_1}$ and ${\grg}_{\alpha_2}$. So $u=0$. 

We are left to prove 1) in the case when $l=k$. If all roots of $\Phi$ have equal lengths, then the above part involving $u_1$ and $u_2$ for $DT\neq A_{2n}$ applies entirely. If $DT$ is $F_4$ or $G_2$, then $H^{\ad}=H^{\sc}$ and so [13, Hauptsatz] applies. So to prove 1) for $l=k$, we can assume $DT\in\{B_n,C_n|n\Ge 2\}$. We can also assume $p=2$, cf. 3.4. As in the mentioned part involving $u_1$ and $u_2$ we get that $[u,\grg_{\beta}]=\{0\}$, provided $\beta\in\Phi$ is long. For any short root $\alpha\in\Phi$ there is a long root $\beta\in\Phi$ such that $\alpha+\beta\in\Phi$; so the component of $u$ in $\grg_{\alpha}$ with respect to (5) is $0$, cf. 3.6 b). For any long root $\beta\in\Phi$ there is a short root $\alpha\in\Phi$ such that $\alpha+\beta\in\Phi$ is also short. So the component of $u$ in $\grg_{\beta}$ is also $0$, as otherwise the component of $AD_H(g)(u)-u$ in $\grg_{\alpha+\beta}$ is non-zero (cf. 3.6 c)); here $g\in\dbG_{a,\alpha}(k_1)$ is an arbitrary non-identity element. So $u\in\Lie(T)$. As $H$ is adjoint, all $\tilde H(|\alpha|)$ groups are $PGL_2$ groups (cf. 3.8). So from (7) and the $DT=A_1$ case we get that $u\in\Lie(T(|\alpha|))$, $\forall\alpha\in\Phi$. As $H$ is adjoint we have $\cap_{\alpha\in\Phi} \Lie(T(|\alpha|)=\{0\}$. So $u=0$. This ends the proof of 1).

To prove 2) we can assume $H$ is s.c. We assume such a $V$ does exist and we show that this leads to a contradiction. We first show that we can assume $k_1=k$. We write $H_{k_1}=\tilde H_1\times_{k_1}\tilde H_2$, where $\tilde H_1$ and $\tilde H_2$ are normal, s.c. subgroups of $H_{k_1}$, the adjoint of $\tilde H_1$ being simple. Let $\tilde V_1:=\Ker(V\otimes_k k_1\to\Lie(\tilde H_2^{\ad}))$. Any trivial $\tilde H_1(k_1)$-submodule of $\tilde V_1$ is contained in $\Lie(Z(H_{k_1}))$, cf. 1). We get $\tilde V_1=\tilde V_1\cap\Lie(\tilde H_1)+\tilde V_1\cap\Lie(Z(H_{k_1}))$. From this and the inclusion $\Lie(\tilde H_1)\subset\tilde V_1+\Lie(Z(H_{k_1}))$ we get $\Lie(\tilde H_1)=\tilde V_1\cap\Lie(\tilde H_1)+\Lie(Z(\tilde H_1))$. So if we have $\Lie(\tilde H_1)=\tilde V_1\cap\Lie(\tilde H_1)$ for any such $\tilde H_1$, then $V\otimes_k k_1=\Lie(H_{k_1})$ and this contradicts the fact that $V$ is a proper $k$-vector subspace of $\Lie(H)$. So we can assume $k_1=k$; so $H=H_1$. If $DT$ is not (resp. is) $D_{2n+2}$, then $\dim_{k}(\Lie(Z(H)))$ is $1$ (resp. is $2$). So the number $\dim_k(\Lie(H)/V)$ is $1$ (resp. is $1$ or $2$). We can assume $\dim_k(\Lie(H)/V)=1$. 

We consider first the case $p\Ge 3$. So either $DT=A_{pn-1}$ or $p=3$ and $DT=E_6$. If $l=k$, then $\dim_{k}(\Lie(H(|\alpha|))\cap V)\Ge 2$ and so $\Lie_k(H(|\alpha|))\cap V=\Lie_k(H(|\alpha|))$ (cf. 3.4), $\forall\alpha\in\Phi$. So $\Lie(H)=\grL_H\subset V$, cf. 3.7 1) and 2). Contradiction. If $l\neq k$, then $l=k_2$. Using similar intersections $V\cap \Lie_k(H_1(i))$ we get (cf. Cases 1, 2 and 3 of 3.9) that to prove $V=\Lie(H)$ we can assume that $DT=A_2$, $p=3$ and we just need to show that $\Lie(H_1(2))\subset V$. The $H_1(2)(k)$-module $\Lie(H)$ is semisimple: it is the direct sum of $\Lie(H_1(2))$ with $\Lie(Z(H))$ and with two $k$-vector spaces of dimension $2$ included in $\grg_{-\alpha_1}\oplus\grg_{-\alpha_2}\oplus\grg_{\alpha_1}\oplus\grg_{\alpha_2}$. So as $V$ is normalized by $H_1(2)(k)$, by reasons of dimensions we get that $V$ contains these two $k$-vector spaces and $\Lie(H_1(2))$. Let $g\in (U_1\setminus H_1(2))(k)$. Let $x\in V\cap (\grg_{-\alpha_1}\oplus\grg_{-\alpha_2}\oplus\grg_{\alpha_1}\oplus\grg_{\alpha_2})\setminus\{0\}$. It is easy to see that we can choose $g$ and $x$ such that the component of $AD_H(g)(x)-x$ in $\Lie(Z(H))$ is non-zero. So as $AD_H(g)(x)-x\in V$ we get $\Lie(Z(H))\subset V$. So $\Lie(H_1(2))\subset V$. So $V=\Lie(H)$. Contradiction.  

We are left with the case $p=2$. The case $DT\in\{A_1,E_7\}$ follows from [13, Hauptsatz]. So we can assume $DT\in\{A_{2n-1},B_n,C_n|n\Ge 2\}\cup\{D_{n+3}|n\in\dbN\}$. As in the previous paragraph we just need to show that $\Lie_k(H_1(i))\cap V=\Lie_k(H_1(i))$, for any $H_1(i)$ as in 3.9. The case $l=k\neq\dbF_2$ is as above. If $l=k=\dbF_2$, then 2) follows from loc. cit. So we can assume $l\neq k$. So $DT$ is $A_{2n-1}$ or $D_{n+3}$ and we are in one of the Cases 1, 4 or 5 of 3.9. In Case 4 all $H_1(i)$'s are $\Res_{k_2/k}$ of $SL_2$ groups, cf. 3.9.1; so as $l=k_2\neq\dbF_2$ we get $\Lie_k(H_1(i))\cap V=\Lie_k(H_1(i))$. The same applies to Cases 1 and 5, except for the situation when $l=k_3$ and $k_1=k=\dbF_2$; so $DT=D_4$. We now refer to this situation and we use the notations of Case 5 of 3.9. It is enough to show that $V\cap\Lie(H_1(2))=\Lie(H_1(2))$. As $H_1(2)$ is an $SL_2$ group (cf. 3.9.1), we have $\{0,x+h,y+h,x+y\}\subset V\cap\Lie(H_1(2))$, where $x$, $y$ and $h$ are the non-zero elements of $\Lie(H_1(2))$ such that $x\in\grg_{\alpha_1+\alpha_2+\alpha_3+\alpha_4}$, $y\in\grg_{-\alpha_1-\alpha_2-\alpha_3-\alpha_4}$ and $h\in\Lie(T)$. If $g\in\dbG_{a,\alpha_2}(l)\cap H(k)$ is a non-identity element, then $u:=AD_H(g)(x+y)-(x+y)$ is the non-zero element of $\grg_{\alpha_1+2\alpha_2+\alpha_3+\alpha_4}\cap\Lie(H)$ (cf. 3.6 c)). Similarly, if $w\in\dbG_{a,-\alpha_2}(l)\cap H(k)$ is the non-identity element, then $x=AD_H(w)(u)-u$. So $V\cap\Lie(H_1(2))$ contains $x$. So $V\cap\Lie(H_1(2))=\Lie(H_1(2))$. This ends the proof of 2).

We prove 3). As $Z(H_1^{\sc})=\mu_2^2$ and $Z(H_1)=\mu_2$ we have $\dim_{k_1}(\Lie(H_1^{\ad})/\grL_{H_1^{\ad}})=2$ and $\dim_{k_1}(\Lie(H_1)/\grL_{H_1})=1$. So $\Lie(Z(H_1))\subset\grL_{H_1}$. So 3) follows from 2) applied to the inverse image $V$ of $V_0$ in $\Lie(H^{\sc})$. We prove 4). Either $DT=A_{p^2n-1}$ or $p=2$ and $DT=D_{2n+3}$. We have a direct sum decomposition of $H$-modules
$\Lie(H)=\Lie(Z(H))\oplus [\Lie(H),\Lie(H)]=\Lie(Z(H))\oplus [\Lie(H^{\ad}),\Lie(H^{\ad})]$, cf. [18, 0.13]. So 4) follows from 3.7.1. To prove 5) we can assume $l=k$. We have $\dim_k(\Lie(H)/\grL_H)\Le 2$ and the equality can take place only for $p=2$, cf. loc. cit. As $\dim_k(H)\ge 3$ and as for $p=2$ the group $H$ is not an $SL_2$ group, the representation of $H$ on $\Lie(H)/\grL_H$ is trivial. So $H(k)$ acts trivially on $\Lie(H)/\grL_H$. This ends the proof. 

\medskip\noindent
{\bf 3.11. On exceptional ideals.} Until 3.12 we assume that $H$ is adjoint and that either $p=2$ and $DT\in \{B_n,C_n|n\in\dbN\}\cup\{F_4\}$ or $p=3$ and $DT=G_2$. So $l=k_1$. It is known that $H_1$ is the extension to $k_1$ of a split, adjoint group $H_0$ over $k$, cf. [27, Vol. III, p. 410]. Let $I$ (resp. $I^{\sc}$) be the ideal of $\Lie(H_1)$ (resp. $\Lie(H_1^{\sc})$) generated by $\grg_{\alpha}$, with $\alpha\in\Phi$ short. We have $I=I_0\otimes_k k_1$ (resp. $I^{\sc}=I_0^{\sc}\otimes_k k_1$), with $I_0$ (resp. $I_0^{\sc}$) as the similarly defined ideal of $\Lie(H_0)$ (resp. $\Lie(H_0^{\sc})$). As $I^{\sc}$ and $I$ are normalized by $H_1$ we have $I=\im(I^{\sc}\to \Lie(H_1))$. So $I$ is an a.s. $H_1(k_1)$-module, cf. [13, Hauptsatz]. Similarly, $I_0$ is an a.s. $H_0(k_0)$-module. Let $I_k$ (resp. $I_k^{\sc}$) be $I$ (resp. $I^{\sc}$) but viewed as a $k$-ideal or as a $k$-vector subspace of $\Lie(H)$. 

\medskip\noindent
{\bf 3.11.1. Proposition.} {\it The only simple $H(k)$-submodule of $\Lie(H)$ is $I_k$.}

\medskip
\proof 
Let $I_{H_0}\otimes_k k_1:H_0(k)\to GL(I_0\otimes_k k_1)(k)$ be the representation defined by the restriction of $AD_H$ to $H_0(k)\vartriangleleft H_0(k_1)=H_1(k_1)=H(k)$. It is semisimple, being a direct sum of $[k_1:k]$ copies of the irreducible representation $I_{H_0}:H_0(k)\to GL(I_0)(k)$ defined by $AD_{H_0}$. As $I_0$ is an a.s. $H_0(k_0)$-module, any non-trivial subrepresentation of $I_{H_0}\otimes_k k_1$ is of the form $IRR=I_0\otimes_k W$, with $W$ a non-trivial $k$-vector subspace of $k_1$. We choose $B_1$ and $T_1$ such that $T_1$ is the extension to $k_1$ of a maximal torus $T_0$ of $H_0$. So $T_1(k_1)=T_0(k_1)$ is included in $H_0(k_1)$. So $W$ is stable under multiplications with elements of $k_1^*$ which belong to the image of the homomorphism $T_1(k_1)\to\dbG_m(k_1)$ associated to a character $T_1\to\dbG_m$ defining a short root of $\Phi$. So $W$ is $k_1^*$-stable and so $W=k_1$. So $IRR=I_0\otimes_k k_1$. So $I_k$ is a simple $H(k)$-module. From loc. cit. and (4) we get that any simple $H(k)$-submodule of $\grL_H\otimes_k k_1$ is contained in $I_k\otimes_k k_1$. So $I_k$ is the only simple $H(k)$-submodule of $\grL_H$. Any other simple $H(k)$-submodule of $\Lie(H)$ would be trivial (cf. 3.10 5)) and so does not exist (cf. 3.10 1)). This ends the proof. 

\medskip\noindent
{\bf 3.11.2. Theorem.} {\it {\bf 1)} If $DT$ is $B_n$ or $C_n$ we assume $n$ is odd and at least $3$. Then the $H(k)$-module $\grL_H/I_k$ is simple and non-trivial. Moreover, the short exact sequence of $H(k)$-modules $0\to I_k\to\grL_H\to\grL_H/I_k\to 0$ does not split.

\smallskip
{\bf 2)} If $DT=C_n$ with $n$ odd and at least $3$, then there are no 1 dimensional $k$-vector subspaces of $\Lie(H)/I_k$ normalized by $H(k)$.

\smallskip
{\bf 3)} Let $H_1^{\sc\text{long}}$ be the analogue of $H_1^{\text{long}}$ of 3.9 but for $H_1^{\sc}$ instead of $H_1$. We assume $DT=C_n$ (resp. $DT=B_n$). If $n=1$ we also assume $k\neq\dbF_2$. Then $[\Lie(H^{\sc}),\Lie(H^{\sc})]$ (resp. $I_k+\Lie_k(Z(H_1^{\sc\text{long}}))$) is the only maximal $H(k)$-submodule of $\Lie(H^{\sc})$.}

\medskip
\proof
The first part of 1) is argued as in 3.11.1 for $I_k$. The second part of 1) follows from 3.11.1. To prove 2) we can assume $k=k_1$. So $H=H_1$ and $I=I_k$. The group $H^{\text{long}}:=H_1^{\text{long}}$ is the quotient of the group $SL(H)=SL_2^n$ of the proof of 3.7 by $Z(H^{\sc})$. Let $T_{\sc}$ be the inverse image of $T$ in $H^{\sc}$. We have $\dim_k(\Lie(H^{\sc})/I^{\sc})=2n$, cf. [13, Hauptsatz]. This implies $\Lie(T_{\sc})\subset I^{\sc}$ and so $\Lie(H)/I$ (resp. $\grL_H/I$) is $2n+1$ (resp. $2n$) dimensional. So $\Lie(H)/I$ is identified naturally with a Lie subalgebra of $\Lie((H^{\text{long}})^{\ad})$ in such a way that $\grL_H/I$ corresponds to $\grL_{({H^{\text{long}}})^{\ad}}$. But $\Lie((H^{\text{long}})^{\ad})$ has no 1 dimensional $k$-vector subspace normalized by the subgroup $H^{\text{long}}(k)$ of $H(k)$, cf. 3.10 1). So 2) holds.

We now prove 3). We deal just with the $DT=C_n$ case as the $DT=B_n$ case is entirely the same. The $H(k)$-module $\Lie(H^{\sc})/[\Lie(H^{\sc}),\Lie(H^{\sc})]$ is simple. Argument: the case $n=1$ follows from 3.7.1 and the case $n>1$ is argued in the same way 3.11.1 was, the only difference being that $T_1\to\dbG_m$ defines this time a long root of $\Phi$. Let $M$ be an $H(k)$-submodule of $\Lie(H^{\sc})$ not contained in $[\Lie(H^{\sc}),\Lie(H^{\sc})]$. It is enough to show that $M=\Lie(H^{\sc})$. We show that the assumption $M\neq \Lie(H^{\sc})$ leads to a contradiction. Let $\tilde M$ be a  maximal $H_1(k_1)$-submodule of $\Lie(H^{\sc}_{k_1})=\Lie(H^{\sc})\otimes_k k_1$ containing $M\otimes_k k_1$. From (4) and loc. cit. we get that $[\Lie(H^{\sc}),\Lie(H^{\sc})]\otimes_k k_1\subset\tilde M$. As $M$ surjects onto $\Lie(H^{\sc})/[\Lie(H^{\sc}),\Lie(H^{\sc})]$, $\tilde M$ surjects onto $(\Lie(H^{\sc})/[\Lie(H^{\sc}),\Lie(H^{\sc})])\otimes_k k_1$. So $\tilde M=\Lie(H^{\sc}_{k_1})$. Contradiction. So $M=\Lie(H^{\sc})$. This ends the proof.

\medskip
From 3.7.1 and 3.11.2 3) we get:

\medskip\noindent
{\bf 3.12. Corollary.} {\it In general, the $H^{\sc}(k)$-module $\grL_{H^{\ad}}$ has a unique  maximal submodule $\grI_{H^{\ad}}$. Moreover, the $H^{\sc}(k)$-module $\grL_{H^{\ad}}/\grI_{H^{\ad}}$ is non-trivial.}

\bigskip
\noindent
{\boldsectionfont \S 4. The $p$-adic context}

\bigskip
Let $k=\dbF_{q}$ and $s\in\dbN$. Let $G$ be a reductive group scheme over $W(k)$. Let $\dbA^1$ be the affine line over $B(k):=W(k)[{1\over p}]$. The group $G(B(k))$  is endowed with the coarsest topology making all maps $G(B(k))\to\dbA^1(B(k))=B(k)$ induced by morphisms $G_{B(k)}\to\dbA^1$ to be continuous. Let $K$ be a closed subgroup of $G(W(k))$.

\medskip\noindent
{\bf 4.1. Problem.} {\it Find practical conditions on $G$, $K$ and $p$ which imply $K=G(W(k))$.}

\medskip
This Problem was first considered in the context of $SL_n$ groups over $\dbZ_p$ in [24, Ch. IV]. After mentioning two general properties (see 4.1.1 and 4.1.2), in 4.1.3 and 4.1.4 we identify a key subpart (question) of this Problem. In 4.2 we list some simple properties of $G$. In 4.3 we include an approach of inductive nature which allows to get information on $\im(K\to G(W_2(k)))$ from information on semisimple subgroups of $G$. In 4.4 to 4.7 we assume $G$ is semisimple. Section 4.5 solves question 4.1.4 for adjoint and s.c. groups $G$. The computations needed for this are gathered in 4.4. In 4.6 we include a supplement to 4.5 related to the exceptional ideals. In 4.7 we prove 1.3. 

We now point out that in general we do need some conditions on $G$, $K$ and $p$. 

\medskip\noindent
{\bf 4.1.1. Example.} We assume there is an isogeny $f:G^1\to G$ over $W(k)$ of order $p^n$. So $\Ker(f)$ is a finite, flat subgroup of $Z(G^1)$ of multiplicative type (see [27, Vol. III, 4.1.7 of p. 173]). Being of order $p^n$ it is connected. So the Lie homomorphism $\Lie(G^1_k)\to\Lie(G_k)$ has non-trivial kernel and so it is not an epimorphism. Moreover, the group $\Ker(f)(\kbar)$ is trivial and so the homomorphism $f(k):G^1(k)\to G(k)$ is an isomorphism. So the image of the homomorphism $f(W(k)):G^1(W(k))\to G(W(k))$ is a proper subgroup of $G(W(k))$ surjecting onto $G(k)$. 

\medskip
We identify $\Ker(G(W_{s+1}(k))\to G(W_s(k)))$ and $\Lie_{\dbF_p}(G_k)$ as $G(k)$-modules. Let 
$$L_{s+1}:=\im(K\to G(W_{s+1}(k)))\cap\Ker(G(W_{s+1}(k))\to G(W_s(k))).$$
\noindent
{\bf 4.1.2. Lemma.} {\it If $p=2$ (resp. $p>2$) let $m:=3$ (resp. $m:=2$). If $K$ surjects onto $G(W_m(k))$, then $K=G(W(k))$.}

\medskip
\proof
This is just the generalization of [24, Exc. 1 p. IV-27]. We show that in the abstract context of $G$ we can appeal to matrix computations as in loc. cit. Let $\rho_{B(k)}:G_{B(k)}\hookrightarrow GL(W)$ be a faithful representation, with $W$ a finite dimensional $B(k)$-vector space. So $\rho_{B(k)}$ extends to a representation $\rho:G\to GL(L)$, with $L$ a $W(k)$-lattice of $W$ (cf. [19, 10.4 of Part I]). The morphism $\rho$ is a closed embedding, cf. [30, c) of Proposition 3.1.2.1]. So each element $g\in\Ker(G(W(k))\to G(W_s(k)))$ mod $p^{s+1}$ is congruent to $1_L+p^sx$, where $x\in\Lie(G)\subset\End(L)$. So if $s\Ge 3$ (resp. $s\Ge 2$), then the image of $g$ in $G(W_{s+1}(k))$ is the $p$-th power of any element of  $G(W_{s+1}(k))$ whose image in $G(W_{s}(k))$ is the same as the one of $1_L+p^{s-1}x$ mod $p^s$. So by induction on $s\Ge 3$ (resp. on $s\Ge 2$) we get that $L_{s+1}=\Ker(G(W_{s+1}(k))\to G(W_{s}(k)))$ and so that $K$ surjects onto $G(W_{s+1}(k))$. As $s$ is arbitrary and $K$ is compact we get $K=G(W(k))$. This ends the proof.

\medskip
If $p=2$, then 4.1.2 does not hold in general if we replace $m=3$ by $m=2$. Example: $G=\dbG_m$ and $q=2$. See 4.7.1 below for more examples. 

\medskip\noindent
{\bf 4.1.3. The class $\gamma_G$.} We view $\Lie_{\dbF_p}(G_{k})$ as a $G(k)$-module via $AD_{G_k}$. Let 
$$\gamma_G\in H^2(G(k),\Lie_{\dbF_p}(G_{k}))$$ 
be the class defining the standard short exact sequence 
$$0\to\Lie_{\dbF_p}(G_k)\to G(W_2(k))\to G(k)\to 0.\leqno (9)$$
Lemma 4.1.2 points out that for the study of 4.1 we need to understand what possibilities for $L_2$ (resp. for $L_2$ and $L_3$) we have if $p>2$ (resp. if $p=2$). The study of $L_2$ and so of 4.1 is intimately interrelated to the following question. 

\medskip\noindent
{\bf 4.1.4. Question.} {\it When does the short exact sequence (9) have a section (i.e. when is $\gamma_G=0$)?}

\medskip\noindent
{\bf 4.2. Proposition.} {\it {\bf 1)} Any maximal torus $T_k$ of a Borel subgroup $B_k$ of $G_k$ lifts to a maximal torus $T$ of $G$. Moreover, there is a unique Borel subgroup $B$ of $G$ containing $T$ and having $B_k$ as its special fibre.

{\bf 2)} The map $G\to G_k$ establishes a bijection between isomorphism classes of reductive (resp. semisimple) group schemes over $W(k)$ and isomorphism classes of reductive (resp. semisimple) group schemes over $k$. 

{\bf 3)} If $G$ is an adjoint (resp. a s.c. semisimple) group, then it is a product of Weil restrictions of a.s. adjoint groups (resp. of s.c. semisimple groups having a.s. adjoints) over Witt rings of finite field extensions of $k$. So if $G^{\ad}$ is simple, then $G$ is of isotypic Dynkin type.}

\medskip
\proof
As $W(k)$ is complete, the first part of 1) follows from [27, Vol. II, p. 47--48]: by induction on $s\in\dbN$ we lift $T_k$ to a maximal torus $T_{W_s(k)}$ of $G_{W_s(k)}$. To prove the second part of 1) we can assume that $T$ is split and so the existence and the uniqueness of $B$ follow from [27, Vol. III, 5.5.1 and 5.5.5 of Exp. XXII]. See [27, Vol. III, Prop. 1.21 of p. 336--337] for 2). Part 3) follows from 2.3.2 and 2).

\medskip\noindent
{\bf 4.3. An inductive approach.} Until end let $T$ and $B$ be as in 4.2 1). Let $G_{0k}$ be a semisimple subgroup of $G_k$ normalized by $T_k$. Let $T_{0k}$ (resp. $T_{00k}$) be the maximal subtorus of $T_k$ which is a torus of $G_{0k}$ (resp. which centralizes $G_{0k}$). It lifts uniquely to a subtorus $T_0$ (resp. $T_{00}$) of $T$, cf. [27, Vol. II, p. 47--48]. 

\medskip\noindent
{\bf 4.3.1. Lemma.} {\it The group $G_{0k}$ lifts to at most one closed, semisimple subgroup $G_0$ of $G$ normalized by $T$.}

\medskip
\proof
Let $l$ be the smallest field extension of $k$ such that $T_l$ is split. So $T_{W(l)}$ is also split. We consider the Weyl decomposition 
$$\Lie(G_{W(l)})=\Lie(T_{W(l)})\bigoplus\oplus_{\alpha\in\Phi} \grg_{\alpha}\leqno (10)$$ 
of $\Lie(G_{W(l)})$ with respect to $T_{W(l)}$. Warning: whenever $T$ is split (i.e. $k=l$) we drop the left index $W(l)$. Let $\Phi_0:=\{\alpha\in\Phi|\grg_{\alpha}\otimes_{W(l)} l\subset\Lie(G_{0l})\}$. It suffices to prove the Lemma under the extra assumptions that $T$ is split and that $T_{00}$ is trivial (otherwise we replace $G$ by the derived subgroup of the centralizer of $T_{00}$ in $G$). So $G$ is semisimple and the kernel of the natural homomorphism from $T_k$ into the identity component $Aut^0(G_{0k})$ of the group scheme of automorphisms of $G_{0k}$ is finite. We have $Aut^0(G_{0k})=G_{0k}^{\ad}$, cf. [27, Vol. III, Prop. 2.15 of p. 343]. So $T_k$ and $G_{0k}$ have equal ranks. More precisely, $T_k$ is the maximal torus of $G_{0k}$ whose image in $Aut^0(G_{0k})$ is the same as of $T_k$. So $T_0=T$. As $G_{0k}$ is a semisimple group having $T_k$ as a maximal torus, $\Phi_0$ is symmetric and invariant under the natural action of the group $\text{Aut}_{W(k)}(W(l))$ on $\Phi$. If $G_0$ exists, then it is generated by $T$ and by the $\dbG_a$ subgroups of $G$ normalized by $T$ and whose Lie algebras are $\grg_{\alpha}$, with $\alpha\in\Phi_0$. So as $\Phi_0$ is determined uniquely by $G_{0k}$, the Lemma follows.

\medskip\noindent
{\bf 4.3.2. Proposition.} {\it There is a closed, semisimple subgroup $G_0$ of $G$ lifting $G_{0k}$ iff $\Phi_0$ is closed under addition, i.e. we have $\Phi\cap\{\alpha+\beta|\alpha,\beta\in\Phi_0\}\subset\Phi_0$. So $G_0$ exists if for any simple factor $H$ of $G^{\ad}_k$ none of the three conditions of 3.6 b) holds.}

\medskip
\proof
We can assume that $T$ is split and $T_{00}$ is trivial. So $T$ is a maximal torus of $G_0$. If $G_0$ exists, then the direct sum $\Lie(T)\bigoplus\oplus_{\alpha\in\Phi_0} \grg_{\alpha}$ is a Lie subalgebra of $\Lie(G)$. From Chevalley rule in characteristic 0 (see [27, Vol. III, 6.5 of p. 322]) we get that $\Phi_0$ is closed under addition. We assume now that $\Phi_0$ is closed under addition. So $G_0$ exists and is a reductive group scheme, cf. [27, Vol. III, 5.3.4, 5.4.7 and 5.10.1 of Exp. XXII]. So $G_0$ is semisimple as $G_{0k}$ is so. We assume now that for any simple factor $H$ of $G^{\ad}_k$ none of the three conditions of 3.6 b) holds. So as $\Lie(G_{0k})$ is a Lie subalgebra of $\Lie(G_k)$, from the equality part of 3.6 b) we get that $\Phi_0$ is closed under addition. So $G_0$ exists. This ends the proof.

\medskip
We assume $G_0$ exists. So $\Phi_0=-\Phi_0$ is closed under addition. Let $TG_0$ be the closed subgroup of $G$ generated by $T$ and $G_0$. It is a reductive group scheme, cf. loc. cit.

\medskip\noindent
{\bf 4.3.3. Lemma.} {\it We have a unique direct sum decomposition $\Lie(G)=\Lie(TG_0)\oplus V_0$
of $TG_0$-modules. The resulting direct sum decomposition $\Lie_{\dbF_p}(G_k)=\Lie_{\dbF_p}({TG_0}_k)\oplus V_0/pV_0$ is $G_0(k)$-invariant.}

\medskip
\proof
To check the first part we can assume that $T$ is split and $G$ is semisimple. As $\Phi_0$ is closed under addition, $V_0:=\oplus_{\alpha\in\Phi\setminus\Phi_0}\grg_{\alpha}$ is a $\Lie(TG_0)$-module and so a $TG_0$-module. But $V_0$ is the unique supplement of $\Lie(TG_0)$ in $\Lie(G)$ normalized by $T$, cf. (10). So the first part holds. The second part follows from the first part. This ends the proof.

\medskip\noindent
{\bf 4.3.4. Key Lemma.} {\it If $\gamma_{G_0^{\ad}}\neq 0$, then $\gamma_G\neq 0$.}

\medskip
\proof
It is enough to show that if $\gamma_G=0$, then $\gamma_{G_0^{\ad}}=0$. Let $\gamma_{G_0,G}^0\in H^2(G_0(k),\Lie_{\dbF_p}(G_k))$ be the restriction of $\gamma_G$ via the monomorphism $G_0(k)\hookrightarrow G(k)$. The  component of $\gamma_{G_0,G}^0$ in $H^2(G_0(k),\Lie_{\dbF_p}({TG_0}_k))$ with respect to the decomposition of the second part of 4.3.3, is the $0$ class. So $TG_0(W_2(k))$ has a subgroup $\scrS_0$ mapping isomorphically into $G_0(k)$. The order $o$ of $Z(G_0)(k)$ is prime to $p$. So the image of $\scrS_0$ in $(TG_0)^{\ad}(W_2(k))=G_0^{\ad}(W_2(k))$ is a subgroup mapping isomorphically into $\im(G_0(k)\to G_0^{\ad}(k))$. The index of this last group in $G_0^{\ad}(k)$ is $o$. So as $\Lie_{\dbF_p}(G_{0k}^{\ad})$ is a $p$-group, we get $\gamma_{G_0^{\ad}}=0$. This ends the proof.

\medskip\noindent
{\bf 4.3.5. Lemma.} {\it We assume that $q=2$ and $\gamma_G\neq 0$. Then $\gamma_{G_{W(\dbF_4)}}\neq 0$.}

\medskip
\proof
Let $a\in\dbF_4\setminus k$. We have a direct sum decomposition $\Lie_{\dbF_2}(G_{\dbF_4})=\Lie(G_k)\oplus a\Lie(G_k)$ of $G(k)$-modules. So the restriction of $\gamma_{G_{W(\dbF_4)}}$ via the monomorphism $G(k)\hookrightarrow G(\dbF_4)$ is the direct sum of $\gamma_G$ and of a class in $H^2(G(k),a\Lie(G_k))$. So $\gamma_{G_{W(\dbF_4)}}\neq 0$.

\medskip\noindent
{\bf 4.4. Computations.} 
Until end we assume $G$ is semisimple. For $y\in W_2(k)$ let $\bar y\in k$  be its reduction mod $p$. For $\alpha\in k^*$ let $t_{\alpha}\in W_2(k)$ be the reduction mod $p^2$ of the Teichm\"uller lift of $\alpha$. We start answering 4.1.4 for adjoint and s.c. groups. In 4.4.1 we include a general Proposition which among other things allows us to assume that $q\le 4$. In 4.4.2 (resp. 4.4.3 to 4.4.9) we deal with four (resp. seven) special cases corresponding to $q\le 4$ and s.c. (resp. a.s. adjoint) groups $G$ of small rank. All of them are minimal cases in the sense that 4.4.1 does not apply to them and moreover for any closed, semisimple subgroup $G_0$ of $G$ normalized by $T$ and different from $G$, we have $\gamma_{G_0^{\ad}}=0$. Sections 4.3 and 4.4 allow us to list in 4.5 all cases when $\gamma_G=0$ and $G$ is adjoint or s.c. In other words, all other minimal such cases are ``handled" by 4.4.5, 4.3.5 and the available literature (see proof of 4.5). 

The identity (resp. zero) $n\times n$ matrix is denoted as $I_n$ (resp. as $0_n$). For $i$, $j\in\{1,...,n\}$ let $E_{ij}$ be the $n\times n$ matrix whose all entries are 0 except the $ij$ entry which is $1$. We list matrices by rows in increasing order (the first row, then the second row, etc.). 

\medskip\noindent
{\bf 4.4.1. Proposition.} 
{\it We have $\gamma_G\neq 0$ if any one of the following three conditions hold:

\medskip
-- $q\ge 5$;

\smallskip
-- $q=3$ and $G^{\ad}_k$ has a non-split simple factor which is not a $PGU_3$ group;

\smallskip
-- $q\in\{2,4\}$ and $G^{\ad}_k$ has simple factors which do not split over $\dbF_4$ and are not the $\Res_{W(\dbF_4)/W(k)}$ of a $PGU_3$ group.} 
\finishproclaim

\proof 
We can assume $G^{\ad}$ is adjoint and simple. Our hypotheses imply that $G$ has a closed subgroup $G_0$ normalized by $T$ and isogenous to the Weil restriction of an $SL_2$ group over the Witt ring of a finite field with at least 5 elements. Argument: we can assume $G$ is a.s. (cf. 4.2 3)) and we can deal just with special fibres (cf. 4.2 2) and 4.3.2); so the statement follows from  3.9 and 4.3.2. So based on 4.3.4, if $p\Ge 3$ (resp. $p=2$) it suffices to show that $\gamma_G\neq 0$ for $q\Ge 5$ and $G$ an $SL_2$ (resp. a $PGL_2$) group. If $p\Ge 5$ and $y\in M_2(W_2(k))$, then $(I_2+E_{12}+py)^p=I_2+pE_{12}\neq I_2$. So $G(W_2(k))$ has no element of order $p$ specializing to a non-identity $k$-valued point of a $\dbG_a$ subgroup of $G$. So $\gamma_G\neq 0$ if $p\Ge 5$.

Let now $p=3$; so $r\Ge 2$. For $\alpha\in k^*$ let $X_{\alpha}:=I_2+t_{\alpha}E_{12}+3x_{\alpha}E_{21}+3Y_{\alpha}\in M_2(W(k))$, 
with $x_{\alpha}\in W_2(k)$ and with $Y_{\alpha}\in M_2(W_2(k))$ such that its $21$ entry is $0$. Let $\beta\in k^*$. The matrix equations $X_{\alpha}^3=X_{\beta}^3=I_2$ and $X_{\alpha}X_{\beta}=X_{\beta}X_{\alpha}$ inside $M_2(W_2(k))$, get translated into equations with coefficients in $k$ involving the reductions mod $p$ of $x_{\alpha}$, $x_{\beta}$ and of the entries of $Y_{\alpha}$ and $Y_{\beta}$. We need just few of these equations. Identifying the $12$ (resp. $11$) entries of the equations $X_{\alpha}^3=X_{\beta}^3=I_2$ (resp. $X_{\alpha}X_{\beta}=X_{\beta}X_{\alpha}$) we get $1+{\alpha}\bar x_{\alpha}=1+{\beta}\bar x_{\beta}=0$ (resp. ${\alpha}\bar x_{\beta}={\beta}\bar x_{\alpha}$). So $\alpha^2=\beta^2$. But as $r\Ge 2$ there are elements $\alpha$, $\beta\in k^*$ such that $\alpha^2\neq\beta^2$. So $\gamma_G\neq 0$ if $G$ is an $SL_2$ group, $p=3$ and $r\Ge 2$. 

Let now $p=2$; so $r\Ge 3$. As the homomorphism $GL_2(W_2(k))\to PGL_2(W_2(k))$ is an epimorphism, we can use again matrix computations. For $\delta\in k$ let 
$$X_{\delta}:=I_2+t_{\delta}E_{12}+2Y_{\delta},$$
with $Y_{\delta}\in M_2(W_2(k))$. Let $(a_{\delta},b_{\delta})$ and $(c_{\delta},d_{\delta})$ be the rows of $Y_{\delta}$. As $r\Ge 3$ we can choose $\alpha$, $\beta\in k\setminus\{0,1\}$ such that $1+\alpha+\beta\neq 0$ and $\alpha\neq\beta$. The conditions that the images of $X_1$, $X_{\alpha}$, $X_{\beta}$ in $PGL_2(W_2(k))$ are of order 2 and commute with each other imply that we have the power 2 equations $X_1^2=(1+2v_1)I_2$, $X_{\alpha}^2=(1+2v_{\alpha})I_2$ and $X_{\beta}^2=(1+2v_{\beta})I_2$ and the commuting equations $X_1X_{\alpha}=(1+2v_{1\alpha})X_{\alpha}X_1$, $X_{1}X_{\beta}=(1+2v_{1\beta})X_{\beta}X_1$ and $X_{\alpha}X_{\beta}=(1+2v_{\alpha\beta})X_{\beta}X_{\alpha}$, where $v_1$, $v_{\alpha}$, $v_{\beta}$, $v_{1\alpha}$, $v_{1\beta}$ and $v_{\alpha\beta}$ belong to $W_2(k)$. Among the equations we get are the following ones: 
$$\bar a_1+\bar d_1=\bar a_{\alpha}+\bar d_{\alpha}=\bar a_{\beta}+\bar d_{\beta}=1;\leqno (11)$$ 
$$\bar v_{1\alpha}=\bar v_{1\beta}=\bar v_{\alpha\beta}=1;\leqno (12)$$
$$1+\bar c_{\alpha}+{\alpha}\bar c_1=1+\bar c_{\beta}+{\beta}\bar c_1=1+{\alpha}\bar c_{\beta}+{\beta}\bar c_{\alpha}=0.\leqno (13)$$ 
The equations (11) are obtained by identifying the $12$ entries of the power $2$ equations (we have $X_{\delta}^2=I_2+2t_{\delta}E_{12}+2t_{\delta}[Y_{\delta},E_{12}]$, etc.). The equations (12) are obtained by identifying the $12$ entries of the commuting equations and by inserting (11) into the resulting equations. The equations (13) are obtained by identifying the $11$ entries of the commuting equations and by inserting (12) into the resulting equations. As $1+\alpha+\beta\neq 0$, the system (13) of equations in the variables $\bar c_1$, $\bar c_{\alpha}$ and $\bar c_{\beta}$ has no solution in $k$. So $\gamma_G\neq 0$ if $p=2$, $r\Ge 3$ and $G$ is a $PGL_2$ group. This ends the proof. 

\medskip\noindent
{\bf 4.4.2. Applications.} Let $W_0:=W(\dbF_4)^2$. We refer to the proof of 4.4.1 with $q=4$, $\alpha\in k\setminus\dbF_2$ and without mentioning $\beta$. The equations $X_1^2=X_{\alpha}^2=I_2$ and $X_1X_{\alpha}=X_{\alpha}X_1$ have no solution, cf. the value of $v_{1\alpha}$ in (12). So the class $\gamma_4\in H^2(SL(W_0)(k),\Lie_{\dbF_2}(GL(W_0/pW_0)))$ defined naturally by the inverse image of $SL(W_0)(k)$ in $GL(W_0)(W_2(k))$ is non-zero. So also $\gamma_G\neq 0$ for $G$ an $SL_2$ group over $W(k)$. For most applications 4.3.4 suffices. We now use $\gamma_4$ to get a variant of 4.3.4 in two situations to be used later on.

First we assume that $q=2$ and that $G$ is an $SU_4$ (resp. $SU_5$) group. Let $G_0$ be such that $G_{0W(\dbF_4)}$ is generated by the $\dbG_a$ subgroups of $G_{W(\dbF_4)}$ normalized by $T_{W(\dbF_4)}$ and corresponding as in Case 2 (resp. Case 1) of 3.9 to the roots $\pm\alpha_1$ and $\pm\alpha_3$ (resp. $\pm\alpha_1$ and $\pm\alpha_4$) of $\Phi$. It is the $\Res_{W(\dbF_4)/W(k)}$ of an $SL_2$ group $G_{02}$. The composite monomorphism $G_0(W(k))=G_{02}(W(\dbF_4))\hookrightarrow G(W(k))\hookrightarrow G(W(\dbF_4))$ corresponds to the direct sum $W$ of two copies of the standard rank 2 representation $W_0$ of $G_{02}$ with the trivial rank $0$ (resp. $1$) representation of $G_{02}$ over $W(\dbF_4)$. We show that the assumption $\gamma_G=0$ leads to a contradiction. The image $\gamma_{G_0,G}^1$ of $\gamma_{G_0,G}^0$ in $H^2(G_{02}(\dbF_4),\Lie_{\dbF_2}(GL(W/2W)))$ is also 0. This image is defined by the natural monomorphisms $\Lie(G_k)\hookrightarrow\Lie_{\dbF_2}(G_{\dbF_4})\hookrightarrow\Lie_{\dbF_2}(GL(W/2W))$
of $G_0(k)=G_{02}(\dbF_4)$-modules. But the component of $\gamma_{G_0,G}^1$ in $H^2(G_{02}(\dbF_4),(\Lie_{\dbF_2}(GL(W_0/2W_0)))^4)$ can be identified with four copies of $\gamma_4$ and so it is non-zero. Contradiction. So $\gamma_G\neq 0$. 

Second we assume that $q=4$ and $G$ is an $SU_3$ group. Let $G_0$ be the unique $SL_2$ subgroup of $G$ normalized by $T$. The composite monomorphism $G_0(W(k))\hookrightarrow G(W(\dbF_{16}))=SL_3(W(\dbF_{16}))$ corresponds to the direct sum of $W_0^2$ with the trivial rank $2$ representation of $G_0$ over $W(k)$. So as in the previous paragraph we get that $\gamma_G\neq 0$.

\medskip\noindent
{\bf 4.4.3. On $PGL_3$ over $\dbZ_3$.} 
We assume that $q=3$ and $G$ is a $PGL_3$ group. We show that $\gamma_G\neq 0$. We work inside $M_3(W_2(k))$. Let $X_1:=I_3+E_{12}+3Y_1$ and $X_2:=I_3+E_{23}+3Y_2$, with $Y_1$, $Y_2\in M_3(W_2(k))$. Let $X_3:=X_1X_2X_1^2X_2^2$. We have $X_3=I_3+E_{13}+3Y_3$, with $Y_3\in M_3(W_2(k))$. It is enough to show that we can not choose $Y_1$ and $Y_2$ such that $X_3^3$ is of the form ${\beta}I_3$, with $\beta\in\Ker(\dbG_m(W_2(k))\to\dbG_m(k))$. As $X_3^3=I_3+3E_{13}+3E_{13}Y_3E_{13}$, it is enough to show that the $31$ entry of $X_3$ (and so also of $Y_3$ mod $3$) is $0$.

Let $=_{31}$ be the equivalence relation on $M_3(W_2(k))$ such that two matrices are in relation $=_{31}$ iff their $31$ entries are equal. For $Z\in M_3(W_2(k))$ we have
$$E_{12}Z=_{31} E_{13}Z=_{31} E_{23}Z=_{31} ZE_{12}=_{31} ZE_{23}=_{31} 0_3.\leqno (14)$$ 
So we compute
$$X_3=_{31} (I_3+E_{12}+E_{23}+3E_{12}Y_2+E_{12}E_{23}+3Y_1+3Y_2+3Y_1E_{23})(I_3+2E_{12}+3E_{12}Y_1+6Y_1+3Y_1E_{12})$$
$$(I_3+2E_{23}+3E_{23}Y_2+6Y_2+3Y_2E_{23})=_{31} (I_3+3Y_1+3Y_2+3Y_1E_{23})(I_3+2E_{12}+3E_{12}Y_1+6Y_1+3Y_1E_{12})$$
$$(I_3+3E_{23}Y_2+6Y_2)=_{31} (I_3+3Y_1+3Y_2+3Y_1E_{23})(I_3+3E_{23}Y_2+6Y_2+2E_{12}+6E_{12}E_{23}Y_2+3E_{12}Y_2+$$
$$+3E_{12}Y_1+6Y_1+3Y_1E_{12})=_{31} (I_3+3Y_1+3Y_2+3Y_1E_{23})(I_3+3E_{23}Y_2+6Y_2+6E_{13}Y_2+3E_{12}Y_2+$$
$$+3E_{12}Y_1+6Y_1)=_{31} I_3+3E_{23}Y_2+6Y_2+6E_{13}Y_2+3E_{12}Y_2+3E_{12}Y_1+6Y_1+3Y_1+3Y_2+3Y_1E_{23}.$$  
\noindent
So $X_3=_{31} I_3+9Y_1+9Y_2=_{31} I_3$ and so $\gamma_G\neq 0$.

\medskip\noindent
{\bf 4.4.4. On $PGSp_4$ over $\dbZ_3$.} 
We assume that $q=3$ and $G$ is a $PGSp_4$ group. To show that $\gamma_G\neq 0$ we follow the pattern of 4.4.3. We choose the alternating form $\psi$ on $W(k)^4$ such that for the standard $W(k)$-basis $\{e_1,...,e_4\}$ of $W(k)^4$ and for $i$, $j\in\{1,...,4\}$, $j>i$, we have $\psi(e_i,e_j)=1$ if $j-i=2$ and $\psi(e_i,e_j)=0$ if $j-i\neq 2$. We take $X_1:=I_4+E_{12}-E_{43}+3Y_1$, $X_2:=I_3-E_{14}-E_{23}+3Y_2$, with $Y_1$ and $Y_2$ as in 4.4.3. Defining $X_3$ as in 4.4.3, we have $X_3=I_3+E_{13}+3Y_3$. As in 4.4.3 it is enough to show that the $31$ entry of $X_3$ is $0$. But the computations of 4.4.3 apply, once we remark that similar to (14) for $Z\in M_4(W_2(k))$ and for the equivalence relation $=_{31}$ on $M_4(W_2(k))$ defined as in 4.4.3, we have $(E_{12}-E_{43})Z=_{31} (-E_{14}-E_{23})Z=_{31} Z(E_{12}-E_{43})=_{31} Z(-E_{14}-E_{23})=_{31} 0_4$. 

So $\gamma_G\neq 0$. Warning: here $X_1$, $X_2$ and $X_3$ are ``related" to the roots $\alpha_1$, $\alpha_1+\alpha_2$ and respectively $2\alpha_1+\alpha_2$ of the $C_2$ Dynkin type; the similar computations for the roots $\alpha_1$, $\alpha_2$ and $\alpha_1+\alpha_2$ do not imply that $\gamma_G\neq 0$.

\medskip\noindent
{\bf 4.4.5. On $PGSp_4$ over $\dbZ_2$.} 
We assume that $q=2$ and $G$ is a $PGSp_4$ group. Let $\psi$ be as in 4.4.4 (but with $q=2$). Let $GSp_4:=GSp(W(k)^4,\psi)$. We show that $\gamma_G\neq 0$. The group $GSp_4(k)$ is the symmetric group $S_6$ (see [1]). It can be checked that $G(W(k))$ has a subgroup isomorphic to $A_6$ and so in this case the following computations are more involved. Let $I:=[\Lie(G^{\sc}_k),\Lie(G^{\sc}_k)]$. We view it as a $GSp_4(k)$-module and so also as a normal subgroup of $GSp_4(W_2(k))$. As $GSp_4(k)=G(k)$ and as $\Lie(Z(GSp_4)_k)\subset I$, it is enough to show that the short exact sequence 
$$0\to\Lie({GSp_4}_k)/I\to GSp_4(W_2(k))/I\to GSp_4(k)\to 0\leqno (15)$$ 
does not have a section. We show that the assumption that (15) has a section leads to a contradiction. Let $\scrS$ be a subgroup of $GSp_4(W_2(k))/I$ mapping isomorphically into $GSp_4(k)$. 

If $U\in GSp_4(W_2(k))$ let $\bar U$ be the reduction mod $p$ of $U$ and let $\tilde U$ be the image of $U$ in $GSp_4(W_2(k))/I$. 
Below $x_i$, $y_i$, $z_i$, $w_i\in W_2(k)$, where $i\in\{1,...,5\}$. Let $X\in GSp_4(W_2(k))$ whose rows are $(1+2x_5\,0\,2x_1\,0)$, $(0\,1+2x_5\,0\,1+2x_2+2x_5)$, $(2x_3\,0\,1\,0)$ and $(0\,2x_4\,0\,1+2x_4)$. Let $Y\in GSp_4(W_2(k))$ whose rows are $(1+2y_5\,0\,1+2y_1+2y_5\,0)$, $(0\,1+2y_5\,0\,2y_2)$, $(2y_3\,0\,1+2y_3\,0)$ and $(0\,2y_4\,0\,1)$. Let $Z\in GSp_4(W_2(k))$ whose rows are $(1+2z_5\,0\,2z_1\,2z_1)$, $(-1+2z_5\,1+2z_5\,0\,2z_2)$, $(2z_3\,0\,1\,1)$ and $(2z_4\,2z_4\,0\,1)$. Let $W\in GSp_4(W_2(k))$ whose rows are $(1+2w_5\,0\,2w_1\,1+2w_5)$, $(0\,1+2w_5\,-1+2w_5\,2w_2)$, $(2w_3\,0\,1\,2w_3)$ and $(0\,2w_4\,2w_4\,1)$. The subgroup of $GSp_4(k)$ generated by $\bar X$, $\bar Y$, $\bar Z$ and $\bar W$ is a $2$-Sylow subgroup. 

The subgroup of $\Lie({GSp_4}_k)$ generated by $E_{13}$, $E_{24}$, $E_{31}$, $E_{42}$ and $E_{11}+E_{22}$ is a direct supplement of $I$. So for $U\in\{X,Y,Z,W\}$, $\tilde U$ is the general element of $GSp_4(W_2(k))/I$ lifting $\bar U$. 
We now choose $\tilde X$, $\tilde Y$, $\tilde Z$ and $\tilde W$ to generate a $2$-Sylow subgroup of $\scrS$. So $\tilde X^2$, $\tilde Y^2$, $\tilde Z^2$, $\tilde W^2$, $(\tilde Z\tilde W)^2$, $(\tilde W\tilde X)^2$, $(\tilde X\tilde Y)^2$ and $(\tilde Z\tilde X)^2$ are all identity elements. 

Both $\bar Z$ and $\bar W$ are associated to short roots of the $C_2$ Lie type. So from (6) we get that they belong to the commutator subgroup of $G(k)$. So also $\tilde Z$ and $\tilde W$ belong to the commutator subgroup of $\scrS$ and so to $G^{\sc}(W_2(k))/I$. So $2z_5=2w_5=0$. Looking at the $31$ entry of $Z^2$ (and so of $\tilde Z^2$) we get $2z_4=0$. Looking at the $24$ entry of $W^2$ we get $2w_3=0$. Looking at the $13$ entry of $(ZW)^2$ we get $2w_4=0$. Looking at the $13$ and $24$ entries of $(WX)^2$ and using that $2w_3=2w_4=2w_5=0$, we get that $2x_4=0$ and respectively that $2x_3+2x_4+2x_5=2$. Looking at the $24$ entry of $X^2$ we get $2x_4+2x_5=2$. So $2x_3=0$ and $2x_5=2$. Looking at the $13$ entry of $Y^2$ we get $2y_3+2y_5=2$. Looking at the $13$ entry of $(XY)^2$ we get $2x_3+2x_5+2y_3+2y_5=2$. So $2y_3+2y_5=0$. As we also have $2y_3+2y_5=2$, we reached a contradiction. So (15) does not have a section and so also $\gamma_G\neq 0$.

\medskip\noindent
{\bf 4.4.6. On $PGL_4$ over $\dbZ_2$.} 
We assume that $q=2$ and $G$ is a $PGL_4$ group. We identify $G^{\sc}=SL(W(k)^4)$. We show that the assumption $\gamma_G=0$ leads to a  contradiction. Let $\scrS$ be a subgroup of $G(W_2(k))$ mapping isomorphically into $G(k)$. So $\scrS$ is a simple group, cf. [11, 2.2.7]. The group $G(W_2(k))/\grL_{G_k}$ is isomorphic to $\scrS\times \dbZ/2\dbZ$, cf. 3.10 5). As there is no epimorphism from $\scrS$ onto $\dbZ/2\dbZ$, the images of $\scrS$ and $G^{\sc}(W_2(k))$ in $G(W_2(k))/\grL_{G_k}$ coincide. So the inverse image $\scrS_1$ of $\scrS$ in $G^{\sc}(W_2(k))$ surjects onto $\scrS$. Let $SP:=\{(i,j)|i,j\in\{1,2,3,4\}\,\text{and}\,i\ne j\}$. For $(i,j)\in SP$ let $u_{ij}\in\scrS_1$ be such that its image in $\scrS$ is the reduction mod $p$ of $I_4+E_{ij}$. We have $u_{ij}^2=\pm I_4$. If $s\in\{1,2,3,4\}\setminus\{i,j\}$, then the $ss$ component of $u_{ij}^2$ is 1. So $u_{ij}^2=I_4$. If $(i,j)$, $(s,t)\in SP$ and $v\in\{1,2,3,4\}\setminus\{i,j,s,t\}$, then similarly by identifying the $vv$ components we get $u_{ij}u_{st}=u_{st}u_{ij}$.

Let $X_1:=u_{23}u_{34}u_{23}u_{34}=I_4+E_{24}+2Y_1$ and $X_2:=u_{13}=I_4+E_{13}+2Y_2$, where $Y_1$, $Y_2\in M_2(W_2(k))$. We know that $X_1$, $X_2\in SL_4(W(k)^4)$, $X_1^2=X_2^2=I_4$ and $X_1X_2=X_2X_1$. We take $x_y\in W_2(k)$, where $x\in\{a,b,...,p\}$ and $y\in\{1,2\}$. As $X_1^2=I_4+2E_{24}+2[Y_1,E_{24}]=I_4$, the rows of $X_1$ are $(1+2a_1\,\, 0\,\, 2c_1\,\, 2d_1)$, $(2e_1\,\, 1+2f_1\,\, 2g_1\,\, 1+2h_1)$, $(2i_1\,\, 0\,\, 1+2k_1\,\, 2l_1)$ and $(0\,\, 0\,\, 0\,\, 3+2f_1)$. Similarly, as $X_2^2=I_4$ we get that the rows of $X_2$ are $(1+2a_2\,\, 2b_2\,\, 1+2c_2\,\, 2d_2)$, $(0\,\, 1+2f_2\,\, 2g_2\,\, 2h_2)$, $(0\,\, 0\,\, 3+2a_2\,\, 0)$ and $(0\,\, 2n_2\,\, 2o_2\,\, 1+2p_2)$. Identifying the $24$ entries of $X_2X_1$ and $X_1X_2$ we get $(1+2h_1)(1+2f_2)+2h_2(3+2f_1)=(1+2f_1)2h_2+(1+2h_1)(1+2p_2)$. So $2f_2=2p_2$. As $2f_2=2p_2$ we have $\text{det}(X_2)=(1+2a_2)(3+2a_2)(1+2f_2)(1+2p_2)=3$. So $X_2\notin G^{\sc}(W_2(k))$. Contradiction. So $\gamma_G\neq 0$.

\medskip\noindent
{\bf 4.4.7. On $PGL_3$ over $W(\dbF_4)$.} 
We assume that $q=4$ and $G^{\ad}$ is a $PGL_3$ group. We show that $\gamma_G\neq 0$. As $c(G^{\sc})$ is odd we can assume $G$ is an $SL_3$ group. Let $\alpha$, $\beta\in k^*$, with $\alpha\neq\beta$. Let $X_1:=I_3+t_{\alpha}E_{12}+2Y_1$ and $X_2:=I_3+t_{\beta}E_{12}+2Y_2$, where $Y_1$, $Y_2\in M_3(W_2(k))$. If $X_1^2=X_2^2=I_3$, then $Y_j$ has the rows $(2a_j\,\, 2b_j\,\, 2c_j)$ $(0\,\, 2+2a_j\,\, 0)$ and $(0\,\, 2h_j\,\, 2i_j)$. Here $j\in\{1,2\}$ and all $a_j$, ..., $i_j\in W_2(k)$. So as $\alpha\neq\beta$ we get that the $12$ entries of $X_1X_2$ and $X_2X_1$ are distinct. So $X_1X_2\neq X_2X_1$ and so $\gamma_G\neq 0$. 

\medskip\noindent
{\bf 4.4.8. On $PGU_3$ over $\dbZ_2$.} 
We show that $\gamma_G=0$ if $q=2$ and $G^{\ad}$ is a $PGU_3$ group. We can assume $G$ is s.c. Let $U$ be the unipotent radical of $B$. As $G_{W(\dbF_4)}$ splits, it is the $SL$ group of $M:=W(\dbF_4)^3$ and $T_{W(\dbF_4)}$ splits. So we can choose a $W(\dbF_4)$-basis $\scrB=\{e_1,e_2,e_3\}$ of $M$ such that the $W(\dbF_4)$-spans $<e_1>$ and $<e_1,e_2>$ are normalized by $B_{W(\dbF_4)}$ and the $W(\dbF_4)$-spans $<e_1>$, $<e_2>$ and $<e_3>$ are normalized by $T_{W(\dbF_4)}$. In what follows the matrices of elements of $G(W(\dbF_4))$ are computed with respect to $\scrB$. We can assume that $\scrB$ is such that the automorphism of $U_{W(\dbF_4)}$ defined by the non-identity element $\tau$ of $\Gal(\dbF_4/k)=\Gal(B(\dbF_4)/B(k))$ takes $A\in U(W(\dbF_4))$ whose rows are $(1\,x\,y)$, $(0\,1\,z)$ and $(0\,0\,1)$ into the element $\tau(A)\in U(W(\dbF_4))$ whose rows are $(1\,\tau(z)\,\tau(xz-y))$, $(0\,1\,\tau(x))$ and $(0\,0\,1)$. We have $U(W(k)):=\{A\in U(W(\dbF_4))|\tau(A)=A\}$. As the $G(k)$-submodule $\Lie(G_k)$ of $\Lie_{\dbF_2}(G_{\dbF_4})$ has a direct supplement and as $U(k)$ is a Sylow 2-subgroup of $G(k)$, it is enough to show that the pull back of the standard short exact sequence $0\to\Lie_{\dbF_2}(G_{\dbF_4})\to G(W_2(\dbF_4))\to G(\dbF_4)\to 0$ via the composite monomorphism $U(k)\hookrightarrow G(k)\hookrightarrow G(\dbF_4)$ has a section. Even better, as the $G(k)$-submodule $\Lie_{\dbF_2}(G_{\dbF_4})$ of $\Lie_{\dbF_2}(GL(M/2M))$ has $\Lie_{\dbF_2}(Z(GL(M/2M)))$ as a direct supplement, it is enough to deal with the image of the resulting short exact sequence via the monomorphism $\Lie_{\dbF_2}(G_{\dbF_4})\hookrightarrow\Lie_{\dbF_2}(GL(M/2M))$. 

Let $a\in\dbF_4$ be such that $\dbF_4=\{0,1,a,a+1\}$. Let $\tilde t_a$ and $\tilde t_{a+1}$ be elements of $W_2(\dbF_4)$ lifting $a$ and respectively $a+1$ and such that we have the following two identities 
$$\tilde t_{a+1}+\tilde t_a=\tilde t_a\tilde t_{a+1}=1.\leqno (16)$$ 
Let $X_1$, $X_2$ and $X_3\in GL(M)(W_2(\dbF_4))$ be defined as follows. The rows of $X_1$ (resp. of $X_2$) are $(1\,1\,\tilde t_a)$, $(2\,3+2\tilde t_a\,3)$ and $(0\,0\,1)$ (resp. are $(1\,\tilde t_a\,\tilde t_a)$, $(2+2\tilde t_a\,1\,\tilde t_{a+1})$ and $(0\,0\,3+2\tilde t_a)$). The rows of $X_3$ are $(3\,2\tilde t_a\,3+2\tilde t_a)$, $(0\,3\,0)$ and $(0\,0\,1)$. It is easy to see that $X_3^2=I_3$ and $X_1^2=X_2^2=(X_1X_2)^2=X_3$. We include here just the only two computations which appeal to (16). The $13$ entry of $X_2^2$ is $\tilde t_a+\tilde t_a\tilde t_{a+1}+3\tilde t_a+2(\tilde t_a)^2=1+(2\tilde t_a)^2=3+2\tilde t_a$. The rows of $X_1X_2$ are $(3+2\tilde t_a\,1+\tilde t_a\,2(\tilde t_a)^2+\tilde t_{a+1})$, $(2\tilde t_a\,3\,3+3\tilde t_{a+1})$ and $(0\,0\,3+2\tilde t_a)$. So the $13$ entry of $(X_1X_2)^2$ is $2(3+2\tilde t_a)[2(\tilde t_a)^2+\tilde t_{a+1}]+(1+\tilde t_a)(3+3\tilde t_{a+1})=\tilde t_{a+1}+3+3\tilde t_a+3\tilde t_a\tilde t_{a+1}=3+2\tilde t_a$. 

So the subgroup of $GL(W)(W_2(\dbF_4))$ generated by $X_1$, $X_2$, $X_3$ is a quaternion group of order 8. Its reduction mod $2$ is $U(k)$. So $\gamma_G=0$ for the present case.

\medskip\noindent
{\bf 4.4.9. On $PGU_3$ over $\dbZ_3$.} 
We assume that $q=3$ and $G^{\ad}$ is a $PGU_3$ group. We show that $\gamma_G\neq 0$. Let $M$, $\scrB$, $U$, $\tau$ be as in 4.4.8 but with $q=3$ and with $\dbF_4$ replaced by $\dbF_9$. Let $a\in\dbF_9$ be such that $a^2+1=0$. Let $X_1$, $X_2$, $X_3$, $Y_1$, $Y_2$, $Y_3\in M_3(W(\dbF_9))$ be such that $X_1=I_3+t_aE_{12}+t_{2a}E_{23}+2E_{13}+3Y_1$, $X_2=I_3+t_{a+1}E_{12}+t_{2a+1}E_{23}+E_{13}+3Y_2$ and $X_3=X_1X_2X_1^2X_2^2=I_3+t_{2a}E_{13}+3Y_3$. The reductions mod $3$ of $X_1$, $X_2$ and $X_3$ define elements of $U(k)$. As in 4.4.3 and 4.4.4 we get that the $31$ entry of $Y_3$ mod $3$ is $0$. So $X_3^3$ is not a scalar multiple of $I_3$. So $\gamma_G\neq 0$.

\medskip\noindent
{\bf 4.5. Theorem.} {\it We assume $G$ is adjoint. Then $\gamma_G=0$ (resp. $\gamma_{G^{\sc}}=0$) iff $q\Le 4$ (resp. $q\Le 3$) and $G$ (resp. $G^{\sc}$) is a product of adjoint (resp. s.c.) groups of the following type:

\medskip
{\bf F2} $PGL_2$, $PGL_3$, $\Res_{W(\dbF_4)/\dbZ_2} PGL_2$, $PGU_3$, $PGU_4$ and split of $G_2$ Dynkin type (resp. $SL_3$, $SU_3$ and split of $G_2$ Dynkin type) if $q=2$;

{\bf F3} $PGL_2$ (resp. $SL_2$) if $q=3$;

{\bf F4} $PGL_2$ if $q=4$.}

\medskip
\proof
If $\gamma_{G^{\sc}}=0$, then $\gamma_G=0$. It suffices to prove the Theorem for a simple, adjoint group $G$,  cf. 4.2 3) and 4.3.4. We first assume $\gamma_G=0$. Let $k_1$ be such that $G=\Res_{W(k_1)/W(k)}G_1$, where $G_1$ is an a.s. adjoint group (cf. 4.2 3)). Let $T_1$ (resp. $B_1$) be the maximal torus (resp. Borel subgroup) of $G_1$ such that $\Res_{W(k_1)/W(k)} T_1=T$ (resp. $\Res_{W(k_1)/W(k)} B_1=B$), cf. 2.3.1. Let $G_{01k_1}$ be a closed, semisimple subgroup of $G_{1k_1}$ normalized by $T_{1k_1}$ and obtained as in 3.9 (i.e. working with respect to $B_{1k_1}$). It lifts to a semisimple subgroup $G_{01}$ of $G_1$ normalized by $T_1$, cf. 4.3.2 and 3.9 (2). So $G_0:=\Res_{W(k_1)/W(k)} G_{01}$ is a closed, semisimple subgroup of $G$ (cf. 2.3.1) normalized by $T$.

As $\gamma_G=0$ we have $p\Le 3$, cf. 4.4.1. If $p=3$, then $q=3$ and $G$ is split (cf. 4.4.1 and 4.4.9). We show that $G$ is a $PGL_2$ group. If it is not, then we can take $G_0$ to be split of $A_2$ or $C_2$ Dynkin type (even if $G$ is of $G_2$ Dynkin type). So $\gamma_{G_0^{\ad}}\neq 0$ (cf. 4.4.3 and 4.4.4) and so $\gamma_G\neq 0$ (cf. 4.3.4). We reached a contradiction. So $G$ is a $PGL_2$ groups. But $\gamma_G=0$ and $\gamma_{G^{\sc}}=0$ for such a $G$, cf. [24, p. IV-27-28].

Until end of the proof we assume $p=2$. As $\gamma_G=0$, $k$ is a subfield of $\dbF_4$ and either $G$ splits over $W(\dbF_4)$ or is the $\Res_{W(\dbF_4)/W(k)}$ of a $PGU_3$ group (cf. 4.4.1). 

\medskip
{\bf i)} The group $G_0^{\ad}$ can not be a $PGL_4$ group over $\dbZ_2$, cf. 4.4.6 and 4.3.4. So also $G_0^{\ad}$ can not be the $\Res_{W(\dbF_4)/W(k)}$ of a $PGL_4$ group, cf. 4.3.4 and 4.3.5. So $G$ is not the Weil restriction of a split group of $E_6$ Dynkin type, cf. 3.9 and 4.3.2. Also $G$ is not a $PGU_{n+7}$ group, a $PGL_{n+3}$ group, or of isotypic $D_{n+3}$, $E_7$ or $E_8$ Dynkin type (cf. 3.9 and 4.3.2). 

{\bf ii)} The group $G_0^{\ad}$ can not be a $PGSp_4$ group over $\dbZ_2$, cf. 4.4.5 and 4.3.4. So also $G_0^{\ad}$ can not be the $\Res_{W(\dbF_4)/W(k)}$ of a $PGSp_4$ group, cf. 4.3.4 and 4.3.5. So $G$ is not of isotypic $B_n$ with $n\Ge 2$, $C_n$ with $n\Ge 2$ or $F_4$ Dynkin type (cf. 3.9 and 4.3.2).

{\bf iii)} The group $G_0^{\ad}$ can not be the $\Res_{W(\dbF_4)/\dbZ_2}$ of a $PGL_3$ group, cf. 4.4.7 and 4.3.4. So $G$ is not a $PGU_6$ or a $PGU_7$ group, cf. 3.9 and 4.3.2.

{\bf iv)} The group $G_0$ can not be a $PGU_6$ group over $\dbZ_2$, cf. iii) and 4.3.4. So $G$ is not non-split of $E_6$ Dynkin type, cf. Case 3 of 3.9 and 4.3.2. 

{\bf v)} The groups $G_0$ and $G$ are not the $\Res_{W(\dbF_4)/W(k)}$ of a $PGL_3$ group, cf. 4.4.7 and 4.3.4. So $G$ is also not $\Res_{W(\dbF_4)/W(k)}$ of a split group of $G_2$ Dynkin type, cf. 3.9 and 4.3.2.

{\bf vi)} The group $G$ is not a $PGU_5$ group or the $\Res_{W(\dbF_4)/W(k)}$ of a $PGU_3$ group, cf. 4.4.2 and the fact that $c(G)$ is odd.

\medskip 
So $G$ is among the adjoint groups listed in F2 and F4. If $q=2$ and $G$ is split of $G_2$ Dynkin type, then the epimorphism $G(W(k))\twoheadrightarrow G(k)$ has a right inverse (cf. [12, Sect. 4]).  It is well known that this also holds if $q=2$ (resp. $q=4$) and $G$ is a $PGL_n$ group with $n\in\{2,3\}$ (resp. $n=2$). Let $SO_6^-$ be the semisimple group over $\dbZ_2$ which is an isogeny cover of degree 2 of a $PGU_4$ group. The epimorphism $SO_6^-(\dbZ_2)\twoheadrightarrow SO_6^-(\dbF_2)$ has a right inverse (see [1, p. 26]; $SO_6^-(\dbF_2)$ is the Weyl group $W_{E_6}$). Based on all these and 4.4.8, we get that F2 and F4 list all cases when $p=2$, $G$ is a simple adjoint group and $\gamma_G=0$. 

The passage from adjoint groups to s.c. groups for $p=2$ is easy. It is well known that $\gamma_{G^{\sc}}\neq 0$ if $q\in\{2,4\}$ and $G$ is a $PGL_2$ group (see 4.4.2 for $q=4$). So $\gamma_{G^{\sc}}=0$ iff $q=2$ and $G^{\sc}$ is a product of s.c. semisimple groups having a.s. adjoints of $A_2$ or $G_2$ Dynkin type (cf. 4.4.2 for the exclusion of $SU_4$ groups for $q=2$). This ends the proof.  

\medskip\noindent
{\bf 4.5.1. Remark.} 
We have $\gamma_G=0$ if $q=2$ and $G$ is the quotient of $\Res_{W(\dbF_4)/\dbZ_2} SL_2$ by $\mu_2$ (this can be deduced either by just adapting (11) to (13) or from [1, p. 26] via the standard embedding of $G$ into $SO_6^-$).

\medskip\noindent
{\bf 4.6. A supplement to 4.5.} 
Until 4.7 we assume $G$ is split, s.c. and $G^{\ad}$ is a.s. We also assume that either $p=3$ and $G$ is of $G_2$ Dynkin type or $p=2$ and $G$ is of $B_n$, $C_n$ or $F_4$ Dynkin type. If $G$ is an $SL_2$ group, then we also assume that $q\Ge 8$. Let $I$ be the maximal $G(k)$-submodule of $\Lie_{\dbF_p}(G_k)$, cf. 3.11.2 1) and 3). We refer to (10). Let $G_0$ be the closed, semisimple subgroup of $G$ generated by the $\dbG_{a}$ subgroups of $G$ having $\grg_{\alpha}$'s as their Lie algebras, with $\alpha\in\Phi$ a long root (see 3.9 and 4.3.2). 

\medskip\noindent
{\bf 4.6.1. Lemma.} {\it We have identifications of Lie algebras}
$$\Lie(G_k)/I=\Lie(G_{0k})/\Lie(Z(G_{0k}))=\grL_{G_{0k}^{\ad}}.\leqno (17)$$
\noindent
\proof 
The intersection $I\cap\Lie(G_{0k})$ is fixed by any maximal torus of $G_{0k}$ and so also by $G_{0k}$. So $I\cap\Lie(G_{0k})\subset\Lie(Z(G_{0k}))$, cf. 3.10 1) applied to $G^{\ad}_{0k}$. A simple computation involving $\dim_k(I)$ (see [13, p. 409]) and $\dim_k(G_{0k})$ (see 3.9 for the Lie type of $G_{0k}^{\ad}$) shows that $\dim_k(I\cap\Lie(G_{0k}))\Ge \dim_k(Z(G_{0k}^{\sc}))$. Example: if $G$ is of $C_n$ Dynkin type, then $\dim_k(I\cap\Lie(G_{0k}))\ge\dim_k(T_k)=n=\dim_k(Z(G_{0k}^{\sc}))$. As $\dim_k(Z(G_{0k}^{\sc}))\Ge \dim_k(Z(G_{0k}))$, we get $I\cap\Lie(G_{0k})=\Lie(Z(G_{0k}))$ and $\dim_k(Z(G_{0k}))=\dim_k(Z(G_{0k}^{\sc}))$. So $\Lie(G_{0k})/\Lie(Z(G_{0k}))=\Lie(G_{0k}^{\sc})/\Lie(Z(G_{0k}^{\sc}))=\grL_{G_{0k}^{\ad}}$ (cf. also proof of 3.10 4) for when $p=2$ and $G_{0k}$ is of $A_3$ or $D_{2n+3}$ Lie type). By reasons of dimensions the monomorphism $\Lie(G_{0k})/\Lie(Z(G_{0k}))=\Lie(G_{0k})/I\cap\Lie(G_{0k})\hookrightarrow \Lie(G_k)/I$ is onto. So (17) holds. This ends the proof.

\medskip\noindent
{\bf 4.6.2. Proposition.} {\it We view $I$ as a normal subgroup of $G(W_2(k))$. Then the short exact sequence
$$0\to\Lie_{\dbF_p}(G_k)/I\to G(W_2(k))/I\to G(k)\to 0\leqno (18)$$
does not have a section.}

\medskip
\proof
Its pull back via the monomorphism $G_0(k)\hookrightarrow G(k)$ is the short exact sequence 
$$0\to\Lie_{\dbF_p}(G_{0k})/\Lie_{\dbF_p}(Z(G_{0k}))\to\im(G_0(W_2(k))\to G_0^{\ad}(W_2(k)))\to G_0(k)\to 0,\leqno (19)$$ 
cf. (17). We show that the assumption that (18) has a section leads to a contradiction. So (19) has a section too. So $\gamma_{G_0^{\ad}}=0$. So $q\Le 4$ and $G$ is of isotypic $C_n$ Dynkin type with $n\Ge 2$ (cf. 4.5; see 3.9 for the structure of $G^{\ad}_{0k}$). To reach a contradiction we can assume $q=2$ (to be compared with 4.3.5 and the proof of 3.11.1). Using a standard monomorphism $\dbG_m^{n-2}\times Sp_4\hookrightarrow Sp_{2n}$ over $W(k)$, similar to 4.3.3 and 4.3.4 we argue that we can assume $n=2$. But as (15) does not have a section (see 4.4.5), (18) does not have a section for $q=2$ and $G$ an $Sp_4$ group. Contradiction. This ends the proof.

\medskip\noindent
{\bf 4.7. Proof of 1.3.} We now prove 1.3. Let $G$ and $K$ be as in 1.3. Let $L_2$ and $L_3$ be as before 4.1.2. As $K$ surjects onto $G(k)$, both $L_2$ and $L_3$ are $G(k)$-submodules of $\Lie_{\dbF_p}(G_k)$. Let $K^{\sc}$ be the inverse image of $K$ in $G^{\sc}(W(k))$. We have $\im(K^{\sc}\to G(k))=\im(G^{\sc}(k)\to G(k))$. So $\im(K^{\sc}\to G^{\sc}(k))$ is a normal subgroup of $G^{\sc}(k)$ of index dividing the order of $Z(G^{\sc})(k)$ and so prime to $p$. So as $G^{\sc}(k)$ is generated by elements of order $p$ (cf. [11, 2.2.6 (f)]), $K^{\sc}$ surjects onto $G^{\sc}(k)$. As $g.c.d.(p,c(G))=1$, the isogeny $G^{\sc}\to G$ is \'etale. So $\Ker(G^{\sc}(W(k))\to G^{\sc}(k))=\Ker(G(W(k))\to G(k))$. So if $K^{\sc}=G^{\sc}(W(k))$, then $K=G(W(k))$. So to show that $K=G(W(k))$ we can assume $G=G^{\sc}$. Let $G^0$ be a direct factor of $G$ such that $G^{0\ad}$ is simple.

We first assume that 1.3 e) holds; so $p=2$ and $L_2=\Lie_{\dbF_p}(G_k)$. To prove that $K=G(W(k))$ it suffices to show that $L_3$ is $\Lie_{\dbF_p}(G_k)=\Ker(G(W_3(k))\to G(W_2(k)))$, cf. 4.1.2. We consider a faithful representation $\rho:G\hookrightarrow GL(L)$ as in the proof of 4.1.2. Let $x$, $y\in\Lie(G)$. Let $h_x$ and $h_y\in K$ be such that we have $h_x=1_L+2x+4\tilde x$ and $h_y=1_L+2y+4\tilde y$, where $\tilde x$, $\tilde y\in \End(L)$. We compute
$$h_xh_yh_x^{-1}h_y^{-1}=(1_L+2x+4\tilde x+2y+4xy+4\tilde y)(1_L+2x+4\tilde x)^{-1}(1_L+2y+4\tilde y)^{-1}=\leqno (20)$$
$$=(1_L+2y+4xy-4yx+4\tilde y)(1_L+2y+4\tilde y)^{-1}=1_L+4(xy-yx).$$
So $L_3$ contains $[\Lie_{\dbF_2}(G_k),\Lie_{\dbF_2}(G_k)]$. So $L_3$ contains also $\Lie_{\dbF_2}(G^0_k)$ except when $G^0$ is of isotypic $C_n$ Dynkin type, cf. 3.7 2). Let now $G^0$ be of isotypic $C_n$ Dynkin type. Let $k_1$ be such that $G^0=\Res_{W(k_1)/W(k)} G_1^0$, with $G_1^0$ a s.c. semisimple group over $W(k_1)$ whose adjoint is a.s. (cf. 4.2 3)). The group $G_1^0$ is split. 

Let $G^0_{10}$ be a closed $SL_2$ subgroup of $G_1^0$ normalized by a maximal split torus $T_1^0$ of $G_1^0$. We identify $G_{10}^0(W_3(k_1))$ with a subgroup of $GL_2(W_3(k_1))$. We denote also by $t_{\alpha}$ the reduction mod $8=p^3$ of the Teichm\"uler lift of $\alpha\in k_1\setminus\{0\}$. Let $j_2$ (resp. $j_4$) be the image in $GL(L)(W_3(k))$ through $\rho$ of the element $I_2+2t_{\alpha}E_{12}$ (resp. $I_2+4t_{\alpha}E_{12}$) of $G^0_{10}(W_3(k_1))\vartriangleleft G(W(k))$. For $y\in\End(L/8L)$ we have $(j_2+4y)^2=j_4\neq I_2$.
So $L_3$ contains the Lie algebra of any $\dbG_a$ subgroup of $G_{10k_1}^0$ normalized by $T_{1k_1}^0$. So as the intersection $L_3\cap\Lie_{\dbF_2}(G^0_{10k_1})$ is $G^0_{10}(k_1)$-invariant, it is $\Lie_{\dbF_2}(G^0_{10k_1})$. So $L_3$ contains $\Lie_{\dbF_2}(G^0_{1k_1})=\Lie_{\dbF_2}(G^0_k)$, cf. 3.7 1) and 2). So $L_3=\Lie_{\dbF_p}(G_k)$ and so $K=G(W(k))$.

We now assume that one of the conditions a) to d) of 1.3 holds. To prove that $K=G(W(k))$ it is enough to show that $L_2=\Lie_{\dbF_p}(G_k)$, cf. 4.1.2 and the above part referring to 1.3 e). We first assume $G^{\ad}$ is simple. Either $\gamma_{G^{\ad}}\neq 0$ or $L_2\not\subset\Lie_{\dbF_p}(Z(G_k))$, cf. 4.5 and our hypotheses. So the group $\tilde L_2:=\im(L_2\to\Lie_{\dbF_p}(G_k^{\ad}))$ is non-trivial. If $G$ is (resp. is not) the Weil restriction of a s.c. group having an a.s. simple adjoint and among the s.c. groups of 4.6, then from 4.6.2 (resp. 3.7.1) we get that $\tilde L_2$ is not included in the unique maximal $G(k)$-submodule of $\im(\Lie_{\dbF_p}(G_k)\to\Lie_{\dbF_p}(G_k^{\ad}))$ (see 3.12). So $\tilde L_2=\im(\Lie_{\dbF_p}(G_k)\to\Lie_{\dbF_p}(G_k^{\ad}))$. So $L_2=\Lie_{\dbF_p}(G_k)$, cf. 3.10 2). So $K=G(W(k))$.

We now treat the general case; so $G^{\ad}$ is not simple. We know that $K$ projects onto $G^0(W(k))$ and that $L_2$ projects onto $\Lie_{\dbF_p}(G^0_k)$. If $x\in L_2$ and $g\in G^0(k)$, then $AD_{G_k}(g)(x)-x\in L_2\cap\Lie_{\dbF_p}(G^0_k)$. From this and 3.12 we get that there is a $G(k)$-submodule of $L_2\cap\Lie_{\dbF_p}(G^0_k)$ surjecting onto $\im(\Lie_{\dbF_p}(G^0_k)\to\Lie_{\dbF_p}(G_k^{0\ad}))$. So from 3.10 2) applied to $\Res_{k/\dbF_p} G^0_k$ we get that the direct factor $\Lie_{\dbF_p}(G^0_k)$ of $\Lie_{\dbF_p}(G_k)$ is contained in $L_2$. As $G^0$ is an arbitrary direct factor of $G$ having a simple adjoint, we get that $L_2$ contains $\Lie_{\dbF_p}(G_k)$ and so it is $\Lie_{\dbF_p}(G_k)$. This ends the proof of 1.3. 

\medskip\noindent
{\bf 4.7.1. Remark.} We assume that $p=2$, $G=G^{\ad}$ and $2|c(G)$. Then there are proper, closed subgroups of $G(W(k))$ surjecting onto $G(W_2(k))$. One checks this using a short exact sequence $0\to\Lie_{\dbF_2}(G_k)/[\Lie_{\dbF_2}(G_k),\Lie_{\dbF_2}(G_k)]\to G(W_3(k))/I\to G(W_2(k))\to 0$ and the fact that for any torus $\tilde T$ over $W(k)$ the square homomorphism $\Ker(\tilde T(W(k))\to \tilde T(k))\to \Ker(\tilde T(W_3(k))\to\tilde T(W_2(k)))$ is not surjective.

\bigskip
\noindent
{\it Acknowledgements.} We would like to thank Serre for two e-mail replies; the first one was the starting point of 4.4 and 4.5 which prove all its expectations and the second one led to a better presentation and to precise references in the proof of 4.5. We would also like to thank D. Ulmer and the referee for many useful comments and U of Arizona for good conditions for the writing of the paper.

\bigskip
\noindent
{\bf References}
\bigskip

{\nspace{

\Ref[1]
J. H. Conway, R. T. Curtis, S. P. Norton, R. A. Parker and R. A. Wilson,
\sl Atlas of finite groups,
\rm Oxford Univ. Press, Eynsham, 1985

\Ref[2]
S. Bosch, W. L\"utkebohmert, M. Raynaud,
\sl N\'eron models,
\rm Springer--Verlag, 1990

\Ref[3]
A. Borel,
\sl Properties and linear representations of Chevalley groups,
\rm LNM {\bf 131}, Springer--Verlag, 1970, pp. 1--55

\Ref[4]
A. Borel,
\sl Linear algebraic groups, 
\rm Grad. Text Math. {\bf 126}, Springer--Verlag, 1991

\Ref[5]
N. Bourbaki,
\sl Groupes et alg\`ebres de Lie, 
\rm Chapitre 4--6, Hermann, 1968

\Ref[6]
N. Bourbaki,
\sl Groupes et alg\`ebres de Lie,
\rm Chapitre 7--8, Hermann, 1975

\Ref[7]
F. Bruhat and J. Tits, 
\sl Groupes r\'eductifs sur un corps local,
\rm Publ. Math. I.H.E.S. {\bf 60},  5--184 (1984)

\Ref[8]
A. Borel and J. Tits, 
\sl Homomorphismes \lq\lq abstraits\rq\rq\ de groupes alg\`ebriques simples, 
\rm Ann. of Math. (3) {\bf 97},  499--571 (1973)

\Ref[9]
C. W. Curtis,
\sl Representations of Lie algebras of classical Lie type with applications to linear groups,
\rm J. Math. Mech. {\bf 9},  307--326 (1960)

\Ref[10]
C. W. Curtis,
\sl On projective representations of certain finite groups,
\rm Proc. Am. Math. Soc. {\bf 11},  852--860, A. M. S. 1960

\Ref[11]
D. Gorenstein, R. Lyons and R. Soloman,
\sl The classification of the finite simple groups, Number 3,
\rm Math. Surv. and Monog., Vol. {\bf 40}, A. M. S. 1994

\Ref[12]
B. Gross,
\sl Groups over $\dbZ$,
\rm Inv. Math. {\bf 124},  263--279 (1996)

\Ref[13]
G. Hiss,
\sl Die adjungierten Darstellungen der Chevalley-Gruppen, 
\rm Arch. Math. {\bf 42},  408--416 (1982)

\Ref[14]
G. M. D. Hogeweij,
\sl Almost-classical Lie algebras,
\rm Indag. Math. {\bf 44}, I.  441--452, II.  453--460 (1982) 

\Ref[15]
J. E. Humphreys, 
\sl Introduction to Lie algebras and representation theory,
\rm Grad. Texts Math. {\bf 9}, Springer--Verlag, 1975

\Ref[16]
J. E. Humphreys, 
\sl Linear algebraic groups,
\rm Grad. Texts Math. {\bf 21}, Springer--Verlag, 1975

\Ref[17]
J. E. Humphreys,
\sl Algebraic groups and modular Lie algebras,
\rm Mem. Am. Math. Soc., No. {\bf 71}, A. M. S., 1976

\Ref[18]
J. E. Humphreys, 
\sl Conjugacy classes in semisimple algebraic groups, 
\rm Math. Surv. and Monog., Vol. {\bf 43}, A. M. S., 1995

\Ref[19] 
J. C. Jantzen, 
\sl Representations of algebraic groups, 
\rm Academic Press 1987

\Ref[20]
R. Pink,
\sl Compact subgroups of linear algebraic groups,
\rm J. of Algebra {\bf 206},  438--504 (1998)

\Ref[21]
K. Ribet, 
\sl On $l$-adic representations attached to modular forms,
\rm Inv. Math., Vol. {\bf 28}, 245--275 (1975)

\Ref[22]
K. Ribet, 
\sl Images of semistable Galois representations,
\rm Pac. J. of Math., Vol. {\bf 3} (3),  277--297 (1997)

\Ref[23] 
J. -P. Serre,
\sl Propri\'et\'es galoisiennes des points d'ordre fini des courbes elliptiques,
\rm Inv. Math. {\bf 15},  259--331 (1972)

\Ref[24]
J. -P. Serre,
\sl Abelian $l$-adic representations and elliptic curves,
\rm Addison--Wesley Publ. Co.  1989

\Ref[25]
 J. -P. Serre, 
\sl Galois Cohomology, 
\rm Springer--Verlag 1997

\Ref[26]
J. -P. Serre, 
\rm Collected papers, Vol. IV, Springer--Verlag 2000

\Ref[27]
M. Demazure, A. Grothendieck, et al. 
\sl Schemes en groupes, 
\rm Vol. I to III, LNM {\bf 151--153}, Springer--Verlag 1970

\Ref[28]
R. Steinberg,
\sl Representations of algebraic groups,
\rm Nagoya Math. J. {\bf 22},  33--56 (1963)

\Ref[29] 
J. Tits, 
\sl Classification of algebraic semisimple groups, 
\rm Proc. Sympos. Pure Math., Vol. {\bf 9},  33--62, A. M. S. 1966

\Ref[30]
A. Vasiu, 
\sl Integral canonical models of Shimura varieties of preabelian type, 
\rm Asian J. Math., Vol. {\bf 3} (2),  401--518 (1999)

}}

\enddocument